\documentclass[11pt,leqno]{article}
\usepackage{amsmath}
\usepackage{amssymb}
\usepackage{amscd}
\usepackage[dvips]{graphics}
\usepackage{amsfonts}
\usepackage{amsfonts}

\newcommand{\refm [1]} {(\ref{#1})}

\newcounter{theorem}

\newtheorem{theorem}{Theorem}[section]
\newtheorem{lemma}[theorem]{Lemma}
\newtheorem{definition}[theorem]{Definition}
\newtheorem{corollary}[theorem]{Corollary}

\newtheorem{remark}[theorem]{Remark}

\newtheorem{proposition}[theorem]{Proposition}

\newcommand{\be}{\begin{equation}}
\newcommand{\eu}{\end{equation}}
\newcommand{\ber}{\begin{eqnarray}}
\newcommand{\ena}{\end{eqnarray}}
\newcommand{\nin}{\noindent}
\newcommand{\non}{\nonumber}

\newcommand{\half}{\frac{1}{2}}

\def\bbb#1{{\rm I\mkern-3.5mu #1}} 
\def \RR {\bbb R}
\newcommand{\mn}{\mu_N}

\newcommand{\ds}{\displaystyle}
\newcommand{\Z}{Z_N}

\newcommand{\E}{\mathbf{E}}
\newcommand{\la}{\label}

\newcommand{\om}{\Omega_N}

\newcommand{\PP}{\mathbb{P}}

\newcommand{\yn}{Y^{(j)}}

\newcommand{\wt}{\varphi^{(j)}(t)}
\newcommand{\wz}{\varphi_{1}^{(j)}}

\newcommand{\ii}{I^{(j)}}
\newcommand{\Var}{{\bf Var}}

\newcommand{\nicef}{{\bf {\cal F}}}

\newcommand{\niceb}{{\bf {\cal B}}}

\newcommand{\nicel}{{\bf {\cal L}}}

\def \R {\bbb R}
\def \RR {\bbb R}

\title{ Limit shapes of Gibbs distributions on the set of integer partitions:
The expansive case
 }

\author{{\bf Michael M. Erlihson}
\thanks{E-mail: maerlich@tx.technion.ac.il}\\
Department of Mathematics, Technion-Israel Institute of Technology,\\
Haifa, 32000, Israel.\\
\quad {\bf Boris L. Granovsky}
\thanks{E-mail: mar18aa@techunix.technion.ac.il} \\
Department of Mathematics, Technion-Israel Institute of Technology,\\
Haifa, 32000, Israel.}
\begin{document}
\maketitle \vskip 5cm \nin American Mathematical Society 2000
subject classifications.

\nin Primary-60J27; secondary-60K35, 82C22, 82C26.

\nin Keywords and phrases: Gibbs distributions on the set of
integer partitions, Limit shapes, Random combinatorial structures
and coagulation-fragmentation processes, Local and integral
central limit theorems.
 \newpage
 \pagestyle{myheadings}
 \markboth{ }{\qquad \qquad\qquad\qquad\qquad\qquad\qquad\qquad\qquad Limit
 shapes}
 \begin{abstract}
 We find  limit shapes for  a family of multiplicative
   measures on the set of partitions, induced by
  exponential generating functions with expansive parameters, $a_k\sim Ck^{p-1},\
    k\rightarrow\infty,\ p>0$, where $C$ is a positive constant. The
  measures considered are associated with the generalized Maxwell-Boltzmann
  models
  in statistical mechanics,
  reversible coagulation-fragmentation processes
   and  combinatorial
 structures, known as assemblies. We prove a  central limit theorem for
 fluctuations of a properly
  scaled  partition chosen randomly according to the above measure, from  its
 limit shape. We  demonstrate that when the component size passes beyond the threshold
 value, the independence of numbers of components
 transforms into their conditional independence (given their masses). Among other things, the paper also discusses, in a general setting, the
 interplay between limit shape, threshold and gelation.
 \end{abstract}
\section{Introduction and Summary}

\ \ \ Given a  sequence of probability measures $\{\mu_N,\
N\geq1\}$
  on the sets of unordered partitions of integers $N\geq 1,$ the limit shape,
  provided it exists, defines the limiting structure, as $N\to \infty$
  of  properly scaled partitions chosen randomly according to the above sequence of measures. The study of the asymptotic structure of random partitions is stimulated by applications to
 combinatorics, statistical mechanics, stochastic processes, etc. Our paper focuses on
 limit shapes for a class of measures associated with the generalized Maxwell-Boltzmann models in statistical mechanics, reversible
 coagulation-fragmentation processes and
 combinatorial structures called assemblies. In the course of  the asymptotic analysis
 of the above class of measures $\mu_N,$ we reveal some interesting phenomenon that, as we expect, will be also seen in models
 related to other multiplicative measures on the set of partitions.

  We describe now  the context of the present paper. In Section 2 we deal with an arbitrary probability measure on the set
 of partitions. We study here the linkage between  the following three important concepts in statistical mechanics and
 combinatorics: threshold, gelation and limit shape. As a by-product of this study, we
 establish the non-existence of
 limit shapes for some known models. Section 3 gives  a definition of a multiplicative
 measure  and outlines the four main fields of  applications of these measures.
    Section 4 contains the statements of  our  main results that are related to
     multiplicative measures  induced by
  exponential generating functions with expansive parameters, $a_k\sim Ck^{p-1},\
    k\rightarrow\infty,\ p>0$, where $C$ is a positive constant.
    Namely, Theorems 4.1 and 4.3 determine the asymptotics of  component  counts
 of sizes that are of order $N^{\frac{1}{p+1}}$ (=the threshold value), while Theorem 4.6
 accomplishes the same for components of small sizes
 (=$o\left(N^{\frac{1}{p+1}}\right)$). In particular, in Theorem 4.3
  we obtain limit shapes of Young diagrams under the measures considered.
   As a corollary of the above three theorems, we reveal  that
 when the component size  passes beyond the above threshold value, the asymptotic independence of
   the component counts  transforms
 into their conditional independence, specified in Theorem 4.1. Section 5 provides
 proofs, that are based on a far
 reaching generalization of
 Khitchine's probabilistic  method. In the last Section 6 we discuss the  limit shapes
 for our models
  versus
 the ones obtained  by Vershik for the generalized models of Bose-Einstein
 and Fermi-Dirac.
The two latter models correspond to a class of multiplicative
measures induced by Euler type generating functions.

\section{The interplay between Limit shape, Threshold and Gelation}

\ \ \ We shall work with the set $\om=\{\eta\}$ of all unordered
partitions
 $\eta=(n_1,\ldots,n_N): \sum_{k=1}^N kn_k=N,$
 of an integer $N$. Here $n_k=n_k(\eta)$ is the number of summands (=components)
 equal to  $k$  in a partition $\eta\in\om.$
 Each $\eta\in\om$ can be depicted by its Young diagram (see e.g. \cite{A}).
  The boundary of a Young diagram (shortly, Young diagram)  of  $\eta\in\om$ is a non-increasing step function
 $\nu=\nu(\bullet;\eta)$ which is given by
 \be
 \nu (u)=\nu(u;\eta)=\sum_{k= u}^N{n_k(\eta)},\quad  u\ge 0, \quad \eta\in
 \om,
  \la{85}
  \end{equation}
 where we set $n_0(\eta)\equiv 0, \ \eta\in \om.$ By the above definition \refm[85], for any integer $u,$ the decrement $\nu(u-0)-\nu(u+0)$ equals the number (possibly $0$) of components
of size $u,$ $\nu(0)$ gives the total number of components,
whereas  the largest $u$ with $\nu(u)>0$ equals the
size of the largest component.\\
  Obviously, $\int_{0}^{\infty}{\nu(u)du}=N, \ \eta\in\Omega_N$. Let now
   $r=r_N>0$ and for a given partition $\eta\in \Omega_N$
   define the scaled Young diagram
  $\widetilde{\nu}=\widetilde{\nu}(\bullet;\eta)$ with the scaling (scaling factor) $r_N:$
  \be
  \widetilde{\nu}(u)=\frac{r_N}{N}\nu(r_N u),\quad u\ge  0,\quad \eta\in\om.
  \la{86}
  \end{equation}

When a partition $\eta\in\om$ is chosen randomly according to a
given probability measure  on $\om,$ it is natural to ask if there
exists a scaling $r_N=o(N),$ such that the random curve
$\widetilde{\nu}(u)$ converges, as $N\to \infty,$ to some
nonrandom curve. Such a curve, if it  exists, is called the limit
shape (of a random Young diagram). To give a formal definition, we
will need some more notations.
  Let $\nicel=\{l(\cdot)\}$ be the space of nonnegative non-increasing functions $l$ on $[0,\infty)$ with
  $\int_{0}^{\infty}{l(u)du}=1$. Clearly, $\widetilde\nu\in \nicel,$ for all
  $\eta\in \om.$ We supply $\nicel$ with the topology of uniform convergence on compact sets in $[0,\infty).$
  For a given $r_N$ denote by $\rho_{r_N}$
   the mapping  $\eta\rightarrow \tilde\nu$ of $\om$ onto $\nicel.$    Given
   a probability measure $\mu_N$ on $\om,$ the mapping $\rho_{r_N}$ induces the measure
  $\rho_{r_N}\mu_N$ on $\nicel,$ $(\rho_{r_N}\mu_N)(l):=\mu_N(\rho_{r_N}^{-1}l),$  $l\in\nicel.$

  In this section, we refer to $\mu_N$ as an arbitrary probability measure on  $\om.$ The definition of a limit shape
  given below follows the one by Vershik in \cite{V1}.
  \begin{definition}
  A continuous curve $l\in\nicel$ is called the \emph{limit shape} of
   a random Young diagram w.r.t. a sequence of measures $\{\mu_N,\ N\geq 1\}$ on
   $\{\Omega_N,\ N\geq1\}$ (the limit shape of $\mu_N$)
   under the scaling $r_N:r_N=o(N),\ N\to\infty$ if the sequence of measures
   $\{\rho_{r_N}\mu_N, N\ge 1\}$
  on $\nicel$ weakly converges, as $N\to\infty,$ to the delta measure which is concentrated on the curve $l$.
  \end{definition}

 Clearly, the weak convergence of the sequence of  measures $\{\rho_{r_N}\mu_N, N\ge 1\}$ in
Definition 2.1 is equivalent to convergence  in probability of
$\widetilde{\nu}(u), \ u\ge 0$ w.r.t. the sequence $\{\mu_N, N\ge
1\}.$ Namely,
   the continuous function $l(\cdot)\in\nicel$ is the limit shape of the
  measure $\mu_N$  under a scaling $r_N$ if for any $0<a<b<\infty$ and
  $\epsilon>0$,
  \ber
  \lim_{N\rightarrow\infty}{\mn\{\eta\in\om:\sup_{u\in[a,b]}{|\widetilde{\nu}(u)-l(u)|}<\epsilon\}}=1.
  \la{1}
  \ena
\refm[1] expresses the law of large numbers for the random process
$\widetilde{\nu}(u), \ u\ge 0.$

 \nin  Note, that if a measure $\mu_N$ has a limit shape $l\in\nicel $ under
 a scaling $r_N,$ then for any $c>0$ the
 function $cl(cu),\ u\geq0$ is a limit shape of the same measure
 $\mu_N$ under the scaling $cr_N$.

 \nin We will say that a measure has no limit shape if there is no scaling $r_N=o(N),$ that
 provides $\refm[1]$ for some $l\in\nicel $.

\nin  {\bf A sketch of the history of limit shapes.}
  The evolution of shapes of random ensembles of particles, as the number of particles goes
  to infinity, was studied for a
  long time in a variety of applied fields: statistical mechanics
  (the Wulf construction for the formation of crystals,
  see
  \cite{shl}, \cite{S1}), stochastic processes on lattices
  (the  Richardson model, see \cite{dur}),
  biology (growth of colonies) etc etc. A special study was concentrated on limit shapes for
  random structures on the set of partitions, in view of applications to statistical mechanics,
  combinatorics, representation theory and additive number systems. In 1977 two independent teams of researchers,
  Vershik $\&$ Kerov \cite{VK1} and Shepp $\&$ Logan \cite{LS}, derived the limit shape of a Young diagram
  w.r.t. the Plancherel measure. Following this seminal result, Pittel \cite{Pit} found the limit shape of
  Young tableaux w.r.t. a uniform  measure. Since the number of Young tableaux corresponding
  to a given partition (Young diagram)  is known to be equal to the degree of the irreducible representation
  associated with the partition, the above uniform measure, as well as the Plancherel
  measure, is related to the hook formula. It should be mentioned
  that the research in this direction revealed also a  deep linkage to the random matrix
  theory, which is  now a rapidly growing subject (see \cite{OK}).
  Parallel to this line of research, Vershik \cite{V1}
  developed a general theory of limit  shapes for a class of measures
  he called  multiplicative and which are discussed in Section 3. These measures encompass
  a wide scope of models from  statistical mechanics and combinatorics, but do not include the measures associated with the
  hook formula. The  results on
   limit shapes of multiplicative measures obtained by Vershik and
   Yakubovich \cite{V23}- \cite{VFY},\cite{Ya},\cite{VY2}
  during the last decade concern  measures induced by
   Euler type generating functions. (We
 refer to some  details of this  research in the course of the present paper).
Extending these  results, Romik \cite{R} derived limit shapes for
multiplicative measures corresponding to some restricted integer
partitions.
 Note that  that the limit shape of the uniform
 measure on the set of partitions (which is a multiplicative measure)
 was firstly obtained via a heuristic argument, by Temperley \cite{T}.
A comprehensive study of this case was done by Pittel in
\cite{Pit1}.

 \nin In contrast, the
 multiplicative measures $\mu_N$ considered in our paper
 are associated with  exponential
generating function.

 This being  said, it should be stressed that the results on limit shapes
 of multiplicative measures were
 stimulated by the  remarkable  papers of Erd$\ddot{o}$s, Turan,  Szalay as
 well as by other  researchers,  on statistics related to integer partitions
 (for more details see \cite{Frist},\cite{Pit},\cite{Pit1},\cite{V1}). In the course of
research on limit shapes of multiplicative measures,  links of the
subject to various fields of
 mathematics were revealed. In particular,  recently the application of probabilistic methods to the study of
 logical limit laws was implemented in \cite{GS} and \cite{ST}. This link is
 based on fundamental theorems of Compton that were extended and
 deepen by Burris and Bell \cite{BB},\cite{B}. The
 place of limit shapes in this latter field is not yet understood.

  It turns out that the existence of a limit shape is closely related to two other phenomena known in statistical mechanics, which are
  gelation and threshold. We recall

 \begin{definition}
  Let  $q_N(\eta)$ be the size of the largest component in a partition $\eta=(n_1,\ldots,n_N) \in\om,$ i.e.
  $q_N(\eta)=\max\{i:n_i>0\}$.

 \nin (i) We say that a measure $\mu_N$   exhibits gelation if for some $0<\alpha<1$,
 \be
  \lim\sup_{N\rightarrow\infty}{\mu_N\{\eta\in\om:\:q_N(\eta)>\alpha N\}}>0.
 \la{23}
 \end{equation}

  \nin (ii) A sequence $\bar{q}_N=o(N),$ $N\to\infty$ is called a threshold for
  the size of the largest component, under  a measure $\mu_N$ (shortly,
  threshold of a measure
  $\mu_N$),  if
   \be
   \lim_{N\rightarrow\infty}\mu_N\left\{\eta\in\om:\:q_N(\eta)\leq\bar{q}_N\right\}=1,
   \la{144}
  \end{equation}
 \nin while  for any sequence $v_N$, s.t. $v_N= o(\bar{q}_N),\ N\to\infty\ $,
 \be
\lim\sup_{N\rightarrow\infty}\mu_N\left\{\eta\in\om:\:q_N(\eta)\leq
v_N\right\}<1.
    \la{1444}
    \end{equation}
    \end{definition}

  \begin{remark} (i) In physics and chemistry, gelation is viewed as a formation of a gel which is a two-phase system consisting of a solid and a liquid in more solid form than a solution. In combinatorics, equivalent names for gelation  are connectedness of components and appearance of a giant component.  We notice that in \cite{kaz}, in contrast to our definition \refm[23], the definition of gelation (=formation a giant component) requires that
 \be
 \lim_{N\rightarrow\infty}{\mu_N\{\eta\in\om:\:q_N(\eta)>\alpha N\}}=1,
 \la{nj}
 \end{equation}
 for some $\alpha>0.$ In \cite{kaz},  some sufficient conditions for the absence of
 gelation (in the above sense) are given for combinatorial models encompassed into Kolchin's
 generalized allocation scheme (see \cite{Kol}).

(ii) It is common (see e.g. \cite{frgr2}) to seek  a threshold
$\bar{q}_N$, if it exists, in the form
$\bar{q}_N=N^{\bar{\beta}},$ where the exponent
$\bar{\beta}:=\inf\{ \beta: \mu_N(\eta\in \om: q_N(\eta)\le
N^\beta )=1\}.$
 \end{remark}

 \nin We  write $(L),(G), (T)$ to abbreviate the statements: "There exists a limit
 shape/ there exists gelation/ there exists threshold",
 respectively and write $(\bar{\bullet})$ to denote the negations of the above statements.

 \begin{proposition}
 \nin Let $\mu_N$ be a measure on $\om$. Then
 \be
 (\bar G)\Leftrightarrow (T), \quad and \quad (L)\Rightarrow (\bar G).
 \la{1445}
 \end{equation}
 Moreover, if \  $\mu_N$ has a limit shape under a scaling $r_N$, then $\mu_N$ has a
 threshold \
 $\bar{q}_N\ge O(r_N),\
 N\to\infty.$
\end{proposition}

  \nin  {\bf Proof:}

 \nin By the definition of gelation,
 \be
 (\bar G)\Leftrightarrow\lim_{N\rightarrow\infty} {\mu_N\left\{\eta\in\om:q_N(\eta)>\epsilon N\right\}}=0,
 \la{90}
 \end{equation}

 \nin for any $\epsilon>0$. Therefore,
 \be
 (\bar G)\Leftrightarrow\lim_{N\rightarrow\infty}{\mu_N\left\{\eta\in\om:q_N(\eta)\le\epsilon_N N\right\}}=1,
 \end{equation}
 \nin for some sequence $\epsilon_N\to 0,\ N\to \infty.$
 This proves the existence of a threshold when there is no
 gelation.
 For  the proof of the second implication in \refm[1445],
our strategy will be to show that $(L)$ implies that for large $N$
the major part of the  total mass $N$ is partitioned among
component sizes of order $O(r_N)$.
 For this purpose, we employ
  the following argument. For given $d\ge 1, u_1>0,$  let
 $u_0=0,\ u_i=iu_1,\ i=1,\ldots, d$ be $d+1$  equidistant nodes, and  assume that $l\in\nicel$ is a limit shape of $\mu_N$. Then
    \ber
    \frac{1}{N}\sum_{k=u_1r_N}^{u_dr_N-1}{kn_k(\eta)}=\frac{1}{N}\sum_{i=1}^{d-1}\sum_{k=u_ir_N}^{u_{i+1}r_N-1}{kn_k(\eta)}
    \geq\sum_{i=1}^{d-1}{u_i\left(\frac{r_N}{N}\sum_{k=u_ir_N}^{u_{i+1}r_N-1}{n_k(\eta)}\right)},\quad\eta\in\om.
    \la{76}
    \ena
    Next, \refm[1] gives  for all $1\le i\le
    d-1$
    \be \lim_{N\rightarrow\infty}{\mu_N\left\{\eta\in\om:
u_i\frac{r_N}{N}\sum_{k=u_ir_N}^{u_{i+1}r_N-1}{n_k(\eta)}>
u_i\Big(l(u_i)-l(u_{i+1})-\epsilon_i\Big)\right\}}=1,\la{st}
\end{equation}

\nin for arbitrary $\epsilon_i>0,\ 1\le i\le d-1$. Consequently,
substituting $u_i=iu_1, \ 1\le i\le d,$ gives

    \be
    \lim_{N\rightarrow\infty}{\mu_N\left\{\eta\in\om:\frac{1}{N}\sum_{k=u_1r_N}^{u_dr_N-1}{kn_k(\eta)}>\sum_{i=1}^{d-1}
    {\big(l(u_i)-l(u_{i+1})\big)u_i}- u_1\sum_{i=1}^{d-1}
i\epsilon_i\right\}}=1.
    \la{2}
    \end{equation}
    Now we  write
    \ber
    \non &  &\sum_{i=1}^{d-1}{\big(l(u_i)-l(u_{i+1})\big)u_i}=\sum_{i=1}^{d-1}{(u_i-u_{i-1})l(u_i)}-l(u_d)u_{d-1}=
    \\ & &=\sum_{i=1}^{d-1}{(u_{i+1}-u_{i})l(u_i)}-l(u_d)u_{d-1}>
    \int_{u_1}^{u_{d}}{l(t)dt}-l(u_d)u_{d-1},
    \la{72}\ena
  where the last equation is due to the fact that the points $u_i, \ i=0,\ldots,d$
  are equidistant.
     Since $\int_{0}^{\infty}{l(t)dt}=1$ and the function $l$ is
non-increasing and continuous on $[0,\infty)$, we have that for a
sufficiently large $d$ and a sufficiently small $u_1>0$,
    \be
     \int_{u_1}^{u_{d}}{l(t)dt}>1-\frac{\epsilon}{3},\quad
      l(u_d)u_{d-1}<\frac{\epsilon}{3},
    \la{70}
    \end{equation}
     for any  $\epsilon>0$. Now it is left to couple $u_1$ with the above
     $\epsilon_i, \ 1\le i\le d-1$ and $\epsilon,$ by setting, say, $\epsilon_i=i^{-3}$ and
      $0<u_1<\frac{2\epsilon}{3}(\sum_{i=1}^\infty i^{-2})^{-1},$ to conclude from
       \refm[2] that for any $\epsilon>0$
    \be
    \lim_{N\rightarrow\infty}\mu_N\left\{\eta\in\om:\frac{1}{N}\sum_{k=u_1r_N}^{u_dr_N-1}{kn_k(\eta)}>1-\epsilon\right\}=1.
    \la{84}
    \end{equation}
 \nin  Since $r_N=o(N),\ N\to \infty$, this immediately implies $\bar{G}$ and, therefore, $T$. We also derive from \refm[84],
  \be
    \lim_{N\rightarrow\infty}\mu_N\left\{\eta\in\om:q_N(\eta)<
    u_1r_N\right\}=0,
    \la{841}
    \end{equation}
 which completes the proof.\qquad $\blacksquare$

The expression \refm[84] exposes the following meaning of the
scaling factor $r_N$, that is not immediately seen from the law of
large numbers \refm[1].
 When $N$ is large, almost all mass $N$ is partitioned into component sizes in the range $[u_1r_N,u_dr_N],$
 for a sufficiently small $u_1>0$
and a sufficiently large $u_d<\infty.$ Consequently, the boundary
of the major part of the random Young diagram scaled by $r_N$
acquires, as $N\to \infty,$ the shape of a nonrandom curve $l.$

\nin Proposition 2.4 is applied in Corollary 3.1, in the
 next section to prove the
 non-existence of limit shapes for  certain classes of
 measures on the set $\om.$

 \section{Gibbs distributions and multiplicative measures}
 \subsection{Mathematical setting }
\ \ \ The subject of the present paper will be  the following
class of measures $\mu_N$ on the set $\om$ of integer partitions
of $N:$ \be
  \mu_N(\eta)=({c_N})^{-1}\:\frac{a_1^{n_1}
   a_2^{n_2}\ldots a_N^{n_N}} {{n_1}!{n_2}!\ldots{n_N}!},\qquad
  \eta=(n_1,\ldots,n_N)\in \Omega_N.
  \la{79}
  \end{equation}
  Here $a=\{a_k\}$ is a positive function on the set of integers, called a parameter function, and $c_N=c_N(a_1,\ldots,a_N)$ is the partition function of
  the measure $\mu_N:$
 \be
  c_0=1,\quad
  c_N=\sum_{\eta\in\Omega_N}\frac{a_1^{n_1} a_2^{n_2}\ldots  a_N^{n_N}}{{n_1}!{n_2}!\ldots{n_N}!}\,,
  \quad\quad N\geq1.
  \la{80}
 \end{equation}

\nin Following Pitman (see \cite{BP},\cite{Pi}), we call the
measures of the  form  \refm[79],\refm[80] Gibbs distributions.
Gibbs distributions are incorporated into the following general
construction formulated by Vershik \cite{V1}
 who was motivated by
applications to statistical mechanics. Let $\Omega=\bigcup_{N\ge
1}\om$ be the set of all integer partitions $\eta=(n_1,\ldots,
n_k,\ldots)$ and let $s=\{s_k\}_1^\infty$ be a sequence of
positive functions on the set of integers. We associate with a
partition $\eta\in \Omega_N$ the function $F_{N}(\eta)=\prod_{k=
1}^N s_k(n_k)$ and define the measure $\mu_N$ on $\om$ by
$\mu_N(\eta)=(c_N)^{-1}F_N(\eta),$ where
$c_N=\sum_{\eta\in\om}F_N(\eta)$ is the   partition function of
$\mu_N$. Next, the  family of probability measures $\mu^{(x)}$ on
$\Omega$, depending on a parameter $x>0$ is constructed in such a
way that for each $N\ge 1$ the conditional probability of
$\mu^{(x)}$  given  $\om,$ is the aforementioned measure $\mu_N,$
for all $x>0$ from the domain of definition of the function
$\nicef(x)=\sum_{N\ge 1}c_Nx^N$. Explicitly, $\mu^{(x)}$ is given
by
$$\mu^{(x)}(\eta)=x^{N(\eta)}(\nicef(x))^{-1}F_{N(\eta)}(\eta), \ \
\eta\in\Omega,\quad 0<x:\nicef(x)<\infty,
$$
where $N(\eta)=\sum_{k\ge 1}kn_k$ denotes the number $N$ which is
partitioned by an $\eta\in\Omega$ and $\nicef$ is the partition
function of $\mu^{(x)}$, or equivalently the generating function
for the sequence $\{c_N\}_1^\infty.$ It is not hard to verify that
the measure $\mu^{(x)}$ possesses the required conditioning
property which formally reads  as follows: \be
\mu_N(\eta):=(\mu^{(x)}\vert
\om)(\eta)=\frac{\mu^{(x)}(\eta)}{\mu^{(x)}(\om)},\quad
\eta\in\Omega_N, \la{pr}\end{equation} for all $x$ in the domain
of definition of the function $\nicef.$ Moreover, it follows from
the above definitions that the generating function $\nicef$ is
expressed as the  Cauchy product of the generating functions
$\nicef_k(x)=\sum_{r\ge 0}s_k(r)x^r,\
 \ k\ge 1$ for the sequences $\{s_k(r)\}_{r\ge 0}, \quad k\ge 1:$
 \be
\nicef(x)=\prod_{k\ge 1}\nicef_k(x^k).\end{equation}

Consequently, \be \mu^{(x)}(\eta)=x^{N(
\eta)}\frac{\prod_{k=1}^{N(\eta)}s_k(n_k)}{\prod_{k\ge1}{\nicef_k(x^k)}},\quad
\eta\in\Omega,\quad 0<x:\nicef(x)<\infty \la{mux}\end{equation}
and \be \mu_N(\eta)=(c_N)^{-1}\prod_{k=1}^{N}s_k(n_k),\quad
\eta\in\om. \la{pr1}
\end{equation}
In view of the representation \refm[mux], Vershik calls the family
(with respect to $x$) of measures $\mu^{(x)}$ multiplicative.
Following \cite{VFY}, we will  preserve the same name for measures
$\mu_N$ induced  by multiplicative families of measures
$\mu^{(x)},$  via \refm[pr]. The representation \refm[mux] tells
us the important fact that the random counts $n_k,\ k\ge 1$ are
independent with respect to the  probability product measures
$\mu^{(x)}.$

From the above formulae the following relation between the two
families of measures  $\mu^{(x)}$ and $\mu_N$ holds:

\be \mu^{(x)}=\sum_{N\ge 1}x^N(\nicef(x))^{-1}c_N\mu_N, \quad
0<x:\nicef(x)<\infty.\la{mumu}
\end{equation}

\nin This  says that  measures $\mu^{(x)}$ can be viewed as a
poissonization of measures $\mu_N,$ which is  a standard
 way to deduce the equivalence, as $N\to \infty,$ of canonical and microcanonical
 ensembles in statistical mechanics
 (see \cite{BSU},
\cite{VY1} and references therein).

In the above setting, the following three particular forms of the
functions $s_k,\ k\ge 1$ and corresponding to them generating
functions $\nicef$ are of great interest:

{\bf Case 1.} \ber \non s_k(r)&=&\frac{a_k^r}{r!},\ \ r\geq
0,\quad a_k>0, \ k\ge 1,\\  \nicef(x)&=& \exp{(\sum_{k\ge
1}a_kx^k)},\quad 0<x:\sum_{k\ge 1}a_kx^k<\infty. \la{sf1}\ena By
\refm[mux], \be \mu^{(x)}(\eta)= \prod_{k\ge 1}
\Big(\frac{(a_kx^k)^{n_k}}{n_k!}\exp{(-a_kx^k)}\Big),\quad
\eta=(n_1,\ldots, n_k,\ldots)\in\Omega, \la{mx1} \end{equation}
i.e. $\mu^{(x)}$ is the probability product measure induced by the
sequence of Poisson ($a_kx^k, \ k\ge 1$) random variables. It
follows from \refm[pr1] that in this case the measures $\mu_N, \
N\ge 1$ are the Gibbs distributions defined by
   \refm[79],\refm[80]. The associated generating function
   $\nicef$ is called exponential.

{\bf Case 2.}

\ber \non s_k(r)&=&\binom{m_k+r-1}{r}, \quad r\ge 0,\quad m_k\ge
1,\quad k\ge 1,\\  \nicef(x)&=&\prod_{k\ge
1}\frac{1}{(1-x^k)^{m_k}}, \quad 0<x: \sum_{k\ge 1}m_kx^k<\infty.
\la{sf2}\ena

Consequently, \refm[mux], \refm[pr1] give
 \ber \non \mu^{(x)}(\eta)&=& x^{N(\eta)}\prod_{k\ge 1}(1-x^k)^{m_k},
 \quad \eta\in\Omega,\\
 \mu_N(\eta)&=&(c_N)^{-1}\prod_{k=1}^N\binom{m_k+n_k-1}{n_k},\quad
 \eta\in\om,
\ena
 which says that $\mu^{(x)}$ is the probability product
measure induced by the negative binomial random variables
NB$(m_k,x^k),\  k\ge 1, $ with a free parameter $ 0<x<1.$ The
function $\nicef$ is  called the Euler type generating function,
because in the case $m_k\equiv 1, \ k\ge 1$ it conforms to the
standard Euler generating function for integer partitions. Note
that in the latter case $\mu_N$ is the uniform measure on $\om.$

{\bf Case 3.}

\ber \non s_k(r)&=&\binom{m_k}{r}, \quad 0\le r\le m_k,\quad
m_k\ge 1,\quad k\ge 1,\\  \nicef(x)&=&\prod_{k\ge 1}(1+x^k)^{m_k},
\quad 0<x: \sum_{k\ge 1}m_kx^k<\infty. \la{sf3}\ena

\nin In view of \refm[sf3],  the following notation for the sets
of restricted partitions is needed: \ber \non
\Omega(m_1,\ldots,m_k,
\ldots)&=&\{\eta=(n_1,\ldots,n_k,\ldots)\in\Omega:
 0\le n_k\le m_k,\ \ k\ge 1\},\\
 \Omega_N(m_1,\ldots,m_N)&=&\{\eta=(n_1,\ldots,n_N)\in\Omega_N:
 0\le n_k\le m_k, \ \ 1\le k\le N\}.\ena

\nin  We have from \refm[mux], \refm[pr1] \ber \non
\mu^{(x)}(\eta)&=&
x^{N(\eta)}\frac{\prod_{k=1}^{N(\eta)}\binom{m_k}{n_k}}{\prod_{k\ge
1}(1+x^k)^{m_k}}, \quad \eta\in\Omega(m_1,\ldots,m_k,
\ldots)\\
 \mu_N(\eta)&=&(c_N)^{-1}\prod_{k=1}^N\binom{m_k}{n_k},\quad
 \eta\in\Omega_N(m_1,\ldots,m_N).
\ena Thus, $\mu^{(x)}$ is the product measure induced by the
binomial random variables  Bin$(m_k,x^k), \ k\ge 1,$ with a free
parameter $ 0<x<1.$ In particular, if $m_k\equiv 1, \ k\ge 1,$
then $\Omega(1,\ldots,1, \ldots),\ \Omega_N(1,\ldots,1)$ are the
sets of integer partitions with distinct components, i.e. each
$n_k$ is either $0$ or $1$, and $\mu_N$ is the   uniform measure
on $\Omega_N(1,\ldots,1).$ \vskip.5cm We see that in all the
aforementioned three cases the sequence $\{\mu_N, \ N\ge 1\}$ of
multiplicative measures
 is given by a single  parameter function,
which is $a=\{a_k\}$ in the first case and $m=\{m_k\}$ in the two
other cases. It is known that for a great scope (but not all!) of
applied models, the asymptotics of a parameter function can be
described by a power law with the exponent $(p-1)\in \R,$ by which
we mean that the parameter function behaves asymptotically as
$Ck^{p-1},  \ C>0, \ k\to\infty.$ From the point of view of the
asymptotics of the corresponding models, the following three
ranges of $p$ should be distinguished (see \cite{bar}): the
logarithmic case, $ p=0,$ the expansive case, $ p>0,$ and the
convergent case, $p<0.$ It was found (see \cite{ABT} and
\cite{bar}
 respectively) that in the logarithmic  and
the convergent cases the measures $\mu_N$ exhibit gelation.
 In view of this, Proposition 2.4 leads to the following
 \begin{corollary}
In  the logarithmic and the convergent cases, multiplicative
measures have not  limit shapes.\end{corollary}

\nin  Therefore, multiplicative measures may have  limit shapes in
the expansive case only.

 \subsection{Applications }

\ \ \ We  outline the four main fields  of applications of
multiplicative measures.

\begin{itemize}
\item {\bf Coagulation- fragmentation processes (CFP's).} Given an
integer $N$, a CFP
   is  a continuous-time Markov chain
  on the set $\om=\{\eta\}$ of all partitions of $N$.
  Here $N$ codes the total population of indistinguishable particles partitioned
  into $n_j$ groups of
  size $ j,$ $j=1,\ldots,N$.  The possible infinitesimal (in time)
   transitions  are coagulation
  of two groups (=clusters) of sizes $i$ and $j$ into
  one  group of size $i+j$ and fragmentation of one group of size $i+j$
   into two groups of sizes $i$ and $j$.
If the ratio of these rates is of the form \be
\frac{a_{i+j}}{a_ia_j},\quad 2\leq
  i+j\leq N,
  \la{78}
  \end{equation}
 \nin where $a=\{a_j\}$ is a positive function on the set of integers,
 called a parameter function of a CFP, then the process conforms to
 the classic mean-field  reversible model of clustering formulated in 1970s by Kelly and
 Whittle. (For more details and history see \cite{K}, \cite{W},\cite{DGG},
 \cite{frgr2}, \cite{GR},\cite{EG}).
 It is  known (see \cite{K},\cite{DGG}) that, given a parameter function $a$,
 the invariant measure  $\mu_N$ of the corresponding  CFP
 is exactly the Gibbs distribution \refm[79].
In particular, under  the parameter function  $a_k=k^{p-1}, \ k\ge
1$ among the  possible rates of coagulation and fragmentation are
the ones given by the kernels $K(i,j)= (ij)^{p-1}$ and
$F(i,j)=(i+j)^{p-1}$ respectively. The cases $p=1,$ constant
rates, and $p=2$ are common in chemical applications (see
\cite{DGG}, \cite{H}).
 \item{\bf Decomposable  random combinatorial structures.} A
decomposable structure of size $N$ is a union of indecomposable
 components, so that the counts $n_1,\ldots ,n_N$ of
components of sizes $1,\ldots,N, $ respectively, form an integer
partition of $N.$  Given a sequence of  integers $m=\{m_k\}$, it
is assumed that each component of size $k$ belongs to one of the
$m_k$ types. Three classes of decomposable structures: assemblies,
multisets and selections, encompass the whole universe of
classical combinatorial objects. Assemblies are structures
composed of labelled elements, so that $a_k=\frac{m_k}{k!}, \ k\ge
1,$ multisets are formed from unlabelled elements and selections
are multisets with distinct components. Supposing that a structure
     is chosen randomly  from the finite set of a certain
     class of structures with size $N,$
 the random partition ${\bf K}^{(N)}$ of an integer $N$, called component
spectrum, is induced: $${\bf K}^{(N)}=
(K_1^{(N)},\ldots,K_N^{(N)}): \sum_{k=1}^NkK_k^{(N)}=N,\ N\ge 1,$$
where the random variable $K_k^{(N)}$ represents the number of
components of size $k$ in  the random structure. The distribution
of ${\bf K}^{(N)}$ is described by a probability measure on the
set $\om $ of integer partitions of $N.$ A remarkable fact is
(see\cite{ABT}) that the distributions of ${\bf K}^{(N)}$
corresponding to assemblies, multisets and selections are just the
multiplicative measures $\mu_N$ induced by Poisson, negative
binomial and binomial distributions respectively. Based on the
combinatorial context, it is also common to treat distributions
$\mu_N$ of  ${\bf K}^{(N)}$ as the ones generated by the
conditioning relation (see \cite{ABT},\cite{GR}),  a frame which
is close (but not identical!) to Vershik's formalism.

\item{\bf Statistical mechanics.} In this context, the
multiplicative measures
 $\mu^{x}$ and $\mu_N$ are referred to as macrocanonical and
 microcanonical ensembles (of particles) respectively. For more details on this topic, see
\cite{V1}, \cite{GR}. The multiplicative measures induced by
Poisson, negative binomial and binomial distributions, with
constant parameter functions, provide a mathematical setting for
the three classical models of ideal gas, called Maxwell-Boltzmann
(MB), Bose-Einstein (BE) and Fermi-Dirac (FD) statistics,
respectively. In \cite{V1} the multiplicative measures with
parameter functions $Ck^{p-1}, \ C>0,\ \ p\ge1,\ k\ge 1,$ are
called generalized classical statistics.

\item {\bf CFP's on set partitions} (see \cite{BP},\cite{Pi}).
 This
 comparatively new setting  is based on Pitman's study of
combinatorial models of random set partitions, which  developed
from the Ewens sampling formula and Kingman's coalescence
processes. We assume here that in the preceding setup for CFP's,
particles are labelled by $1,\ldots ,N,$ so that the state space
of a CFP is  the set $\Omega_{[N]}=\{\pi_{[N]}\}$ of all
partitions  $\pi_{[N]}$ of the set $[N]=\{1,\ldots,N\}.$ Denoting
$\vert A_j\vert$ the size of a cluster (block) $A_j\subseteqq [N]$
of $\pi_{[N]},$ we assign to each $A_j,$ a weight $m_{\vert
A_j\vert} $ which depicts the number of possible inner states of
$A_j$, the states can be e.g., shapes (in the plane or in space),
colors, energy levels, etc. This says that to a set partition
$\pi_{[n],k}$ into $k$ given clusters $A_1,\ldots, A_k$ correspond
$\prod_{j=1}^k m_{\vert A_j\vert} $ different states of the CFP
considered, so that the total number of states formed by all
partitions of the set $[N]$ into $k$ clusters is equal to \be
\sum_{\pi_{[N],k}\in \Omega_{[N],k}}\prod_{j=1}^k m_{\vert
A_j\vert}:=B_{N,k},\label{Bnk}\end{equation} where $\Omega_{[N],k}
$ is the set of all  partitions of the set $[N]$ into $k$ blocks
and the number $B_{N,k}$ is known as a Bell polynomial in weights
$m_1,\ldots, m_{N-k+1}.$  Next, for a given $k\ge 1$ the uniform
measure $p_{[N],k}$ on the set $\Omega_{[N],k}=\{\pi_{[N],k}\}$ is
defined: \be p_{[N],k}(\pi_{[N],k})=\frac{\prod_{j=1}^k m_{\vert
A_j\vert}}{B_{N,k}},\quad
\pi_{[N],k}\in\Omega_{[N],k}.\label{gb}\end{equation} In a more
general setting which encompasses a variety of models (see
\cite{BP},\cite{Pi}), the weights $m_j$ in \refm[gb] are allowed
to be arbitrary nonnegative numbers.  By the known combinatorial
relation,  the distribution $p_{[N],k}$ on $\Omega_{[N],k}$
induces the  distribution $p_{N,k}$ of cluster sizes $\vert
A_j\vert, \ j=1,\ldots,k $ on the set $\Omega_{N,k}$ of integer
partitions of $N$ into $k$ positive summands:
  \be
  p_{N,k}(\eta)=(B_{N,k})^{-1}\prod_{j=1}^N(\frac{m_j}{j!})^{n_j}\frac{1}{(n_j)!},
   \quad \eta=(n_1,\ldots,n_N)\in \Omega_{N,k},
   \label{pnk}
 \end{equation}
 where the partition function $B_{N,k}$  defined as in \refm[Bnk] can be rewritten in the following form:
 \be
 B_{N,k}=\sum_{\eta\in\Omega_{N,k}}\prod_{j=1}^n(\frac{m_j}{j!})^{k_j}\frac{1}{(k_j)!}.
\end{equation}

Pitman calls the measures $p_{N,k},\ p_{[N],k}$ microcanonical
Gibbs distributions on the sets $\Omega_{N,k}$ and
$\Omega_{[N],k}$ respectively.  The linkage of $p_{N,k}$ to the
Gibbs distributions  $\mu_N$ in \refm[79]  is expressed via the
conditioning relation:
$$p_{N,k}=\mu_N\vert \Omega_{N,k},$$
where $\mu_N$ is the Gibbs distribution  given by \refm[95] with
$a_j=\frac{m_j}{j!}, \ j\ge 1.$

However,  the generic model associated with set partitions of
$[N]$ involves a wealth of specific problems (see \cite{BP})
arising from treating $p_{[N],k},\ k=1,\ldots,N$ as marginal
distributions at any time $t,$ of irreversible time continuous
markov processes of pure fragmentation (or pure coagulation) on
the state space $\Omega_{[N]}.$

Finally note that Gnedin \& Pitman \cite{GnPit} and Kerov
\cite{Ker} studied exchangeable Gibbs distributions related to the
Ewens sampling formula (the case $p=0$).
\end{itemize}

 \section{Statement of   main results}
\ \ \   Our forthcoming asymptotic analysis  is devoted
exclusively to Gibbs
   distributions in the expansive  case,
  i.e. when the parameter function
 $a$ of the measure $\mu_N$  in \refm[79] is of the form:
 \be
\ a_k\sim Ck^{p-1}, \quad C>0,\quad p>0,\quad k\rightarrow\infty.\
  \la{95}
 \end{equation}

\nin We already mentioned in the previous section that this case
conforms to the generalized MB statistics and to reversible CFP's
with certain transition rates. To explain the relevance of the
case considered to random combinatorial structures, recall that in
the combinatorial context, Gibbs distributions with the parameter
function \refm[95] describe assemblies with the number of
indecomposable components of size $k,$ equal $m_k=a_kk!\sim
Ck^{p-1}k!, \ C>0, \ k\to \infty. $ Examples of such structures
with $p=1,2$ were given in \cite{EG}. Generally speaking, we can
think that in an initial assembly having, say, $\tilde{m}_k$
indecomposable components on $k $ labelled items, each one of the
$\tilde{m}_k$ components is distinguished by one of  $b_k$  types.
For example, in permutations,  each of $(k-1)!$ cycles on $k$
labelled items is colored by one of $b_k$ colors; in  forests,
each of $k^{k-2} $ trees on $k$ labelled vertices belongs to one
of $b_k$ species, etc. In this setting, the resulting assembly is
composed of $m_k=b_k\tilde{m}_k$ components of size $k$. Thus,
under $b_k=k^{\alpha},$ where $ \alpha>0$ in the first example and
$\alpha>3/2$ in the second example, the resulting assembly
satisfies \refm[95].\\
  Throughout the rest of the paper we assume that all random variables
 considered are induced by the Gibbs distributions  $\mu_N$ with the
 parameter function \refm[95].\\

We note that Vershik's approach \cite{V1},\cite{VFY}, \cite{VY} to
asymptotic analysis of the generalized BE and FD models is based
on the expression \refm[pr], which allows to replace the measure
$\mu_N$ by the product measure $\mu^{(x)},$ with a properly chosen
$x\to 1.$ In contrast to it, our strategy is the straightforward
study of the asymptotics of $\mu_N,$ as $N\to\infty, $ using
stratifications of the integer $N$ and the total number of
components. However, both approaches are based
 on common ideas rooted in statistical mechanics and the saddle point method.

  We first introduce some more notations. Consider two sets of random variables $\nu_{j}$ and $ K_j, \ j=0,\ldots,
 q\geq1$ that determine for a given  $\eta\in \Omega_N$ stratifications
 of the total number of components $\nu=\nu(\eta)
 =\sum_{k=1}^N n_k(\eta)$ and of the total mass $N=\sum_{k=1}^N kn_k(\eta),$ respectively. Namely,
  for given integers $q\ge 1$ and $0=M_0<M_1<\ldots<M_q<N< M_{q+1}=N+1$ and a given  $\eta\in\om$, we set
  \ber
  \non
  \nu_{0}&=&\nu_0(M_1,\eta):=\sum_{k=M_0}^{M_1-1}{n_k(\eta)},\qquad\nu_{j}=\nu_j(M_j,\eta):= \sum_{k=M_j}^{N}{n_k(\eta)},
   \\ \non K_0&=&K_0(M_1,\eta):=\sum_{k=M_0}^{M_1-1}{kn_k(\eta)}=N-K_1,\qquad
  K_j=K_j(M_j,\eta):=\sum_{k=M_j}^{N}{kn_k(\eta)},\\ && j=1,\ldots,q
 \la{6}
 \ena
 and denote
 \ber
  \overrightarrow{\nu}=\left(\nu_{0},\ldots,\nu_{q}\right),\quad
   \overrightarrow{K}=\left(K_1,\ldots,K_q\right).
 \la{115}
 \ena
\nin We will refer to $M_j,\ j=1,\ldots, q$ as stratification
points.
  The  random vectors $\overrightarrow{\nu},
\overrightarrow{K} $ induce the decompositions
$\overrightarrow{\nu}^{\star}$ and $\overrightarrow{K}^{\star}$
 of\ the random variables  $\nu$ and $N-K_0$ respectively,
 into sums of $q+1$ and $q$ disjoint parts:\ber
  \non \overrightarrow{\nu}^{\star}&=&\left(\nu_0^{\star},\ldots,\nu_q^{\star}\right):=\left(\nu_{0},\nu_{1}-\nu_{2},\ldots,
 \nu_{q-1}-\nu_{q},\nu_{q}\right),\\ \overrightarrow{K}^{\star}&=&\left(K_1^{\star},\ldots,
   K_q^{\star}\right):=\left(K_{1}-K_{2},\ldots, K_{q-1}-K_{q},K_{q}\right).
 \la{xxxx}
 \ena
  We also set $K_0^{\star}=K_0$. Clearly,
  \be
  \ds\nu_j^{\star}=\sum_{k=M_j}^{M_{j+1}-1}n_k, \quad \ds K_j^{\star}=\sum_{k=M_j}^{M_{j+1}-1}kn_k,\ \
  j=0,\ldots,q,\quad\ds\sum_{j=0}^q {K_j^{\star}}=N.
  \la{400}
  \end{equation}
 Formally, $\overrightarrow{\nu}^{\star}$ and $\overrightarrow{K}^{\star}$ can be viewed as the
 $1$-$1$ linear transformations
 $\niceb:\RR^{q+1}\rightarrow \RR^{q+1}$ and $\niceb_1:\RR^{q}\rightarrow \RR^{q}$ of the random vectors
 $\overrightarrow{\nu}$ and $\overrightarrow{K}$ respectively:
 \ber
   \overrightarrow{\nu}^{\star}=\niceb\overrightarrow{\nu},\quad\overrightarrow{K}^{\star}= \niceb_1\overrightarrow{K}.
   \la{301}
   \ena
\vskip .5cm \nin It is known from \cite{frgr2} that in the
expansive case \refm[95],  the measure $\mu_N$ has a threshold
$N^{\frac{1}{p+1}},\ p>0$  and that the major part of the total
mass $N,$ when $N$ is large, is concentrated on components of
sizes of order $N^{\frac{1}{p+1}}$. Therefore, in view of the note
after Proposition 2.4, we take throughout the rest of the paper
 \be
 r_N=h^{-1}N^{\frac{1}{p+1}},\qquad h={\left(C\Gamma(p+1)\right)}^{\frac{1}{p+1}}
 \la{h}
 \end{equation}
 as the scaling factor for the limit shape in question. Here the constant $h$ is
chosen in order to simplify the expression for the limit shape.

To establish the desired limit theorems, it is required to define
proper scalings of
 the two pairs of random vectors  in \refm[115] and  \refm[xxxx].
 We denote by $\hat{\bullet}$  the corresponding scaled quantities.
 The explicit expressions for $\hat{\bullet}$ will be chosen in accordance with
 the asymptotic
 problem considered.

\nin Firstly, in  the  forthcoming Theorems 4.1 and
 4.3 and Corollary 4.5, we will study
 the asymptotics of the random quantities
  in \refm[6],\refm[400] in the vicinity of the   $r_N$.
 In this  case, we consider the stratifications
 induced by  the points  $M_j=[ u_jr_N],\ j=0,\ldots,q,\ M_{q+1}=N+1,$  where
 $0=u_0<u_1<\ldots< u_q$ do not depend on $N$, and take the
 scalings $\hat{\bullet}$ in the form

\ber
 \non \hat{\nu}_j^{\star}&=&\sqrt{f_0(p-1)}\ \ \frac{\nu_j^{\star}-S_j^{\star}}{\sqrt{S_0}},\qquad
\\
 \quad \hat{K}_j^{\star}&=&\sqrt{f_0(p+1)}\ \ \frac{K_j^{\star}-E_j^{\star}}{\sqrt{V_0}},\qquad
 j=0,\ldots,q
 \la{68}
 \ena
and
 \ber
  \non \hat{\nu}_j&=&\sqrt{f_0(p-1)}\ \ \frac{\nu_j-S_j}{\sqrt{S_0}},\\
  \quad \hat{K}_j&=&\sqrt{f_0(p+1)}\ \ \frac{K_j-E_j}{\sqrt{V_0}},\qquad
  j=0,\ldots,q.
 \la{104}
 \ena
 Here the constants
  \be
  f_j(s)=C\int_{u_j}^{u_{j+1}}x^{s}e^{-x}dx,\quad j=0,\ldots,q, \quad s\ge 0,
  \la{f0}
 \end{equation}
  whereas the quantities $S_j,E_j, V_j$ depend on the parameter $\delta=\delta_N=(r_N)^{-1}:$
  \ber
  \non S_0=\sum_{k=M_0}^{M_1-1}{a_ke^{-\delta k}},\quad E_0=\sum_{k=M_0}^{M_1-1}{ka_ke^{-\delta k}},\qquad
  V_0=\sum_{k=M_0}^{M_1-1}{k^2a_ke^{-\delta k}},\\
  \non S_j= \sum_{k=M_j}^{N}{a_ke^{-\delta
  k}}, \qquad E_j=\sum_{k=M_j}^{N}{ka_ke^{-\delta
  k}}, \qquad  V_j=\sum_{k=M_j}^{N}{k^2a_ke^{-\delta k}},
  \ena
  \be
  j=1,\ldots q.
  \la{19}
  \end{equation}
 Finally, the starred  quantities $S^{\star}_j,\ E_j^{\star} $ in the RHS's of \refm[68] are
defined respectively as the transformations $\niceb$ and
$\niceb_1$ of the corresponding non-starred quantities $S_j, \
E_j.$ Also, we will use the following abbreviation
 for the
 conditional probabilities
 $$\rho\left(\hat{\nu}_j^{\star}|\hat{K}_j^{\star}\right):=
 \PP\Big(\hat{\nu}_j^{\star}=\bullet_{1,j}|\hat{K}_j^{\star}=
 \bullet_{2,j}\Big),$$
  when the particular values $\bullet_{1,j}, \bullet_{2,j}$ are
  not important.

  \nin The  theorem below that deals with the asymptotic behaviour of the scaled
  quantities defined by \refm[68], will be a source of our subsequent results on limit shapes.
 \begin{theorem}
 Under $N\to\infty$ and a given $q\ge 1,$

 {\rm(i)}  The random variables $\hat{\nu}_j^{\star}, \ j=0,\ldots,q$ are conditionally independent given
  $\hat{K}_j^{\star},j=0,\ldots,q$, s.t. $\sum_{j=0}^{q}{K_j^{\star}}=N$.
  Moreover,
  \be
 \rho\left(\overrightarrow{\hat{\nu}^{\star}}|\overrightarrow{\hat{K}^{\star}}\right)=\prod_{j=0}^{q}
  {\rho\left(\hat{\nu}_j^{\star}|\hat{K}_j^{\star}\right)}.
  \la{92}
  \end{equation}

 {\rm(ii)} The $(2q+1)-$ dimensional random vector
$(\overrightarrow{\hat{\nu}^{\star}},
\overrightarrow{\hat{K}^{\star}})=(\hat{\nu}_0^{\star},\ldots,\hat{\nu}_q^{\star},\hat{K}_1^{\star},\ldots,
   \hat{K}_q^{\star}) $
 converges weakly to the $(2q+1)-$ dimensional gaussian random vector with zero mean
 and the covariance matrix
 $\{\vartheta_{mk}^{\star}(q),\ m,q=0,\ldots,2q\}$ given  by
 \begin{small}
   \ber
   \vartheta_{mk}^{\star}(q)=
  \begin{cases}
   f_m(p-1)\mathbf{1}_{m=k}-\frac{f_m(p)f_k(p)}{\Gamma(p+2)},             & \text{if}\ \ 0\leq m,k\leq q,\\
   f_{m}(p)\mathbf{1}_{m=k-q}-\frac{f_{m}(p)f_{k-q}(p+1)}{\Gamma(p+2)}, & \text{if}\ \ 0\leq m\leq q,\; q+1\leq k\leq 2q,\\
   f_{k}(p)\mathbf{1}_{k=m-q}-\frac{f_{k}(p)f_{m-q}(p+1)}{\Gamma(p+2)}, & \text{if}\ \ 0\leq k\leq q,\; q+1\leq m\leq 2q,\\
   f_{m-q}(p+1)\mathbf{1}_{m=k}-\frac{f_{m-q}(p+1)f_{k-q}(p+1)}{\Gamma(p+2)},
      & \text{if}\ \ q+1\leq m,k\leq 2q,
  \end{cases}
   \la{covm2}
   \ena
   \end{small}
where $\mathbf{1}_{i=j}$ is the Kroneker symbol.

 {\rm (iii)} Moreover,  moments of the random vector $(\overrightarrow{\hat{\nu}^{\star}}, \overrightarrow{\hat{K}^{\star}})$
 converge to the corresponding moments of the gaussian random vector in  {\rm(ii)}.
 In particular,
 \be
 \vartheta_{mk}^{\star}(q)=\lim_{N\to \infty}
  \begin{cases}
    Cov(\hat{\nu}_m^{\star},\hat{\nu}_k^{\star}), \ \
  & \text{if}\ \  0\le m,k\le q, \\
      Cov(\hat{\nu}_m^{\star},\hat{K}_{k-q}^{\star}), \ \ &\text{if}\ \  0\le m\le
     q, \ \ q+1\le k\le 2q,\\
      Cov(\hat{K}_{m-q}^{\star},\hat{K}_{k-q}^{\star}), &
     \text{if}\ \  q+1\le k,m\le 2q.
  \end{cases}
 \la{lcov}
\end{equation}

 {\rm(iv)} Let the stratifications be induced by the equidistant points
 \ $u_j, \ \ j=0,\ldots,q$ and let $k=m+s, \ s> 0,
 \ 0\le m,k\le 2q.$ Then  the absolute value of the covariance, $|\vartheta_{mk}^{\star}(q)|,$ monotonically decreases in
  $s.$ In particular, if \ $q\to \infty,\ s\to \infty,$ while $m$ is
  fixed, then
 \ber
  |\vartheta_{mk}^{\star}(q)|=
  \begin{cases}
  O\left(u_k^p\exp{(-u_{k})}\right), & \text{if}\ \ 0\le m<k\le q,\\
   O\left(u_{k-q}^{p+1}\exp{(-u_{k-q})}\right),    & \text{if}\ \
   q<k\le
   2q,\quad s-q\to \infty.
  \end{cases}
 \la{161}
 \ena

 \end{theorem}
 \begin{remark}
 The first part of \refm[161] expresses the exponential decay of correlations
  between the scaled counts of components of different
 sizes, as  the "distance" $s$ between  the sizes goes to infinity. Phenomena of such kind are widely known in equilibrium
 statistical mechanics (see e.g.\cite{gri}, \cite{grm} and references therein). \end{remark}

  Let $l\in\nicel$ be a limit shape of a measure $\mu_N$ and
let   $\tilde{\nu}(u)$ be defined as in \refm[86]. We
  call
 \be
 \Delta(u):= \tilde{\nu}(u)-l(u),\quad u\ge 0
 \la{fl}
 \end{equation}
 the random fluctuation of $\mu_N$ from its limit shape at a point $u.$
To state the next theorem, let    \be
  b_r(u)=C\Gamma(r+1,u), \quad  u\ge 0,\quad
  r>-1,
  \la{26}
\eu
 where $\Gamma(r+1,u)=\int_{u}^{\infty}{e^{-x}x^{r}dx}, \quad  u\ge 0,\quad
  r>-1$ is the incomplete Gamma function.

 \begin{theorem}(Limit shape of  $\mu_N$  and the cental limit theorem for
 scaled fluctuations   from the limit
 shape).
  Under $N\to\infty$ and a given $q\ge 1,$

  {\rm(i)} The random vector $(\hat{\nu}_1,\ldots,\hat{\nu_q})$ defined by \refm[104]
   weakly converges to the $q-$ dimensional gaussian random
  vector with zero mean and the covariance matrix $\{e_{mk}\}_{1}^{q}$ given  by
 \ber
  e_{mk}=b_{p-1}(u_s)-\frac{b_{p}(u_m)b_{p}(u_k)}{\Gamma(p+2)},\qquad m,k=1,\ldots,q,
 \la{113}
 \ena
 where $s=max{(k,m)}.$

{\rm(ii)} Setting $l(u)=\frac{\Gamma(p,u)}{\Gamma(p+1)}, \  u\ge
0$ in \refm[fl], the relation between the scaled quantities
$\hat{\nu}_1,\ldots,\hat{\nu_q}$  and the random fluctuations
$\Delta(u_j)$ is given by \be
\hat{\nu}_j\sim\big(C\Gamma(p+1)Nr_N^{-1}\big)^{1/2}\Delta(u_j),\quad
  j=1,\ldots,q, \quad N\to \infty, \la{relat}\end{equation} for all
$\eta\in\om$.

 {\rm(iii)} The measure $\mu_N$ has the  limit shape
 $$l_{p-1}(u)=\frac{\Gamma(p,u)}{\Gamma(p+1)}, \ \ u\ge 0,\ \ p>0,$$ under the scaling
 $ r_N$ given by \refm[h].
\end{theorem}
\begin{remark}
 {\rm (a)} Recall  that the measure $\mu_N$ considered has
 the threshold $ N^{\frac{1}{p+1}}$. By
 {\rm(ii)} of the above theorem and \refm[h], the threshold turns out to be of
 the same order as
  the scaling $r_N$ for the limit shape. In view of this, we
  believe that the following stronger form of the second part of \refm[1445] in  Proposition 2.4
  is valid:
 For a wide class of multiplicative measures, the existence of a threshold $\bar{q}_N$ implies the
 existence of a limit shape under a  scaling $r_N=O(\bar{q}_N), \ N\to \infty.$

\nin {\rm(b)} In \cite{VFY}, an analog  of our cental limit
theorem for fluctuations
 was established
  for the
 uniform measure on the set of unordered partitions with distinct summands
  (=classical FD statistics). Recall that
  this  multiplicative measure is associated with the
   generating function $\prod_{k\ge 1}(1+x^k)^{-1}.$
   The structure of the covariance matrix
 in \cite{VFY} is similar to the one in \refm[113]. In particular, the exponential decay of covariances is also seen
 there.
\end{remark}

 \begin{corollary} (The functional
 central limit theorem for the scaled number of components).
 Let in the above stratification scheme, $q=1$, $u_1=u$ and denote $\hat{\nu}(u)=\hat{\nu}_1.$
 Then, for all $u\ge 0$, the scaled random variable $\hat{\nu}(u)$ weakly converges to $N(0,e_{11}), $ where $e_{11}=e_{11}(u),\ u\ge 0$ is as in \refm[113].
\end{corollary}
Recalling that $\nu(0)$ is equal to the total number of components
in a random partition, the particular case $u=0$ of the above
corollary recovers the central limit theorem for  $\nu(0)$
established in
 \cite{EG}.

  We now turn to the asymptotic behavior of the counts of  components  of sizes
  $o(\bar{q}_N)$ (=small
 sizes in comparison to the threshold), s.t. $ o(\bar{q}_N)\to\infty,\ N\to \infty,$
  where $\bar{q}_N= N^{\frac{1}{p+1}}$. For this problem we   make use of
the  stratification points $M_j= M_j(N)$, \ \ s.t.
 $$
  M_j=o\left(N^{\frac{1}{p+1}}\right)\to \infty, \ N\rightarrow\infty,\ \ \ \
    \ \ \ds\lim_{N\rightarrow\infty}{\frac{M_{j}}{M_{j+1}}}:=\rho_j<1,\quad
  j=1,\ldots,q,
$$
\be M_0=0, \ M_{q+1}=N+1,
 \la{inc}
 \end{equation}
  while the scaled quantities $\hat{\nu}_j^\star,
 \ \hat{K}_j^\star$\ are taken
 in the form slightly different from the one in \refm[68]:

    \ber
   \non\hat{\nu}_j^{\star}&=&\frac{\nu_j^{\star}-S_j^{\star}}{\sqrt{S_j^{\star}}},\qquad
    j=0,\ldots,q,\\
    \quad \hat{K}_j^{\star}&=&\frac{K_j^{\star}-E_j^{\star}}{\sqrt{V_j^{\star}}},\qquad
    j=0,\ldots,q-1.
    \la{sca}
    \ena
 \begin{theorem}(Asymptotic independence and the central limit theorem for the counts of components of
 small sizes).
 Let the stratification points satisfy \refm[inc] and
 $\overrightarrow{\hat{\nu}^{\star}}=
 (\hat{\nu}_0^{\star},\ldots,\hat{\nu}_q^{\star}),$
 $\overrightarrow{\hat{K}^{\star}}=(\hat{K}_1^{\star},\ldots,\hat{K}_q^{\star})$
  be
 given by \refm[sca]. Then, as  $N\rightarrow\infty$,\\
 {\rm(i)}  The coordinates of each one of the two random vectors $\overrightarrow{\hat{\nu}^{\star}}$ and
 $\overrightarrow{\hat{K}^{\star}}$ are  independent random variables.\\
 {\rm(ii)} The random vector $\overrightarrow{\chi}=(\hat{\nu}_0^{\star},\hat{K}_0^{\star},\hat{\nu}_1^{\star},
 \hat{K}_1^{\star}\ldots,\hat{\nu}_{q-1}^{\star}, \hat{K}_{q-1}^{\star},\hat{\nu}_q^{\star})$ weakly converges to
 the $(2q+1)$ -dimensional gaussian vector with zero mean and the
 covariance matrix having a diagonal  block structure, with
$q$ blocks

 \be
  B_j^{-1}:=\left( \begin{array}{cccc}
     1 &\alpha_j  & \\
     \alpha_j & 1
     \end{array} \right),\quad j=0,\ldots,q-1,\qquad
     B_q^{-1}=1-\alpha_q^2,
    \la{cov1}
    \end{equation}

\nin where   \ber \non \alpha_j&=&
  \frac{\sqrt{p(p+2)}}{p+1}\frac{1-\rho_j^{p+1}}{\sqrt{(1-\rho_j^{p+2})(1-\rho_j^p)}}
  =\lim_{N\to \infty}
 Cov(\hat{\nu}_j^{\star},\hat{K}_{j}^{\star}),\quad
  j=0,\ldots,q-1,\\ \alpha_q^2&=&1-\lim_{N\to \infty}Var(\hat{\nu}
  _{q}^{\star})=
  \frac{\Gamma^2(p+1)}{\Gamma(p)\Gamma(p+2)}=\frac{p^2}{p^2+p}.
\la{cov2} \ena
  \end{theorem}

\begin{remark}
 It was proven in \cite{frgr2} that in the model considered,
  the counts of components of fixed sizes (= the random variables
 $n_{k_1},\ldots, n_{k_l}, \   0\le k_1<\ldots <k_l<\infty, \ l>1$)
 are asymptotically independent. Combining this result
 with our Theorems 4.1 and 4.6 says that when the component size passes
 beyond the threshold value $ N^{\frac{1}{p+1}},$
 the asymptotic independence of component counts transforms into their conditional
 independence (given masses). In this respect, the threshold value can be also viewed as the critical value for the independence of component counts
 in the model considered.
\end{remark}

 \section{Proofs}
 \subsection{Khintchine-type representation formula}
\ \ \  As in \cite{EG},\cite{frgr1}, \cite{frgr2} and \cite {GS},
our tool for the asymptotic problems considered  will be the
  probabilistic method by A.Khintchine, that was introduced in 1950's in his book \cite{Kh}. The idea of the method
  is to construct the representation of the quantity of interest via the probability function of a sum of independent
  integer valued random variables depending on a free parameter, and then to implement a local limit theorem. In the course of time,
  this method, sometimes without mentioning Khintchine's name,  was applied to investigation of a large scope of  asymptotic problems
  arising in statistical mechanics and in enumeration combinatorics
  (see e.g. \cite{Kol},\cite{AT},\cite{EG},\cite{frgr1},\cite{frgr2},
  \cite{fr},\cite{GS},\cite{C},\cite{mu} and references therein).
 The first three subsections of the present section contain preparatory asymptotic analysis towards the proof of our main
 results.\\
 Clearly, the distributions of the random vectors $\overrightarrow{\nu}$ and $\overrightarrow{K}$ are completely
 determined by the measure $\mu_N$. Explicitly, for any given vectors $\overrightarrow{L}=(L_0,\ldots,L_q)$ and
 $\overrightarrow{N}=(N_1,\ldots,N_q)$,
  \ber
   \non\PP(\overrightarrow{\nu}&=&\overrightarrow{L})=\sum_{\overrightarrow{N}}{\PP\left(\overrightarrow{\nu}=
   \overrightarrow{L},\overrightarrow{K}=\overrightarrow{N}\right)}\\ &=&\sum_{\overrightarrow{N^{\star}}}
   {\PP\left(\overrightarrow{\nu^{\star}}=\overrightarrow{L^{\star}},\overrightarrow{K^{\star}}=\overrightarrow{N^{\star}}
   \right)}:=\sum_{\overrightarrow{N^{\star}}}{R(\overrightarrow{L^{\star}},\overrightarrow{N^{\star}})},
  \la{300}
  \ena
  where, in accordance with \refm[xxxx],\refm[301], we denoted
  \ber
  \non\overrightarrow{L^{\star}}&=&(L_0^{\star},\ldots,L_q^{\star})=\niceb\overrightarrow{L},\\
  \overrightarrow{N^{\star}}&=&(N_1^{\star},\ldots,N_q^{\star})=\niceb_1\overrightarrow{N}.
  \la{401}
  \ena
  \refm[xxxx] implies that for any vectors $\overrightarrow{L^{\star}}$ and $\overrightarrow{N^{\star}}$,
  \be
  R(\overrightarrow{L^{\star}},\overrightarrow{N^{\star}})=\PP\left(\overrightarrow{\nu}=\overrightarrow{L},
  \overrightarrow{K}=\overrightarrow{N}\right).
  \la{302}
  \end{equation}
  Consequently, it follows from \refm[79],\refm[80] that
  \be
    R(\overrightarrow{L^{\star}},\overrightarrow{N^{\star}})=c_N^{-1}\sum_{\eta\in\om}{\left(\ds{\prod_{k=1}^{N}{\frac{a_k^{n_k}}{{n_k}!}}}
  \mathbf{1}_{A}\right)},
  \la{14}
  \end{equation}
  where
  \be
  \la{15}
  A=A(\overrightarrow{L^{\star}},\overrightarrow{N^{\star}})=\{\eta\in\om:\overrightarrow{\nu^{\star}}(\eta)
  =\overrightarrow{L^{\star}},\ \overrightarrow{K^{\star}}(\eta)=\overrightarrow{N^{\star}}\}.
  \end{equation}
  Our first goal will be to derive the  Khintchine type representation for the
  probability $R(\overrightarrow{L^{\star}},\overrightarrow{N^{\star}})$.
 Setting $a_0=0,$ we construct  the array of
discrete random variables $\beta_l^{(j)}$  defined  by
  \be
  \PP(\beta_l^{(j)}=k)=\frac{a_k e^{-\delta_jk}}{S_j^{\star}(\delta_j)},\quad
  k=M_j,\ldots,M_{j+1}-1,\quad j=0,\ldots,q,\quad l\ge 1,
  \la{8}
  \end{equation}
  where the stratification points  $M_0,\ldots,M_{q+1}$ are as in \refm[6],
  $\delta_0>0,\ldots,\delta_q>0$ are free parameters and
  \be
  S_j^{\star}=S_j^{\star}(\delta_j)=\sum_{k=M_j}^{M_{j+1}-1}{a_ke^{-\delta_jk}},\quad
   j=0,\ldots,q,
  \la{10}
   \end{equation}
in accordance with \refm[19].  We assume that for a given $0\le
j\le q,$ the random variables $\beta_l^{(j)}, \ l\ge 1$ are i.i.d.
and for a given $l\ge 1,$ the random variables $\beta_l^{(j)}, \
0\le j\le q$ are independent.

  Distributions of the type \refm[8] are widely used in the
asymptotic analysis related to combinatorial structures,
\cite{Kol},\cite{ABT},\cite{mu}. These distributions  firstly
appeared in the papers of Goncharov (1944) (for references see
\cite{Kol}) and in the monograph \cite{Kh} by Khintchine (1950).

  \begin{lemma}(Khintchine's type representation for the
  probabilities
  $R(\overrightarrow{L^{\star}},\overrightarrow{N^{\star}}))$.
  Denote
  \be
  Y^{(j)}=\sum_{l=1}^{L_j^{\star}}{\beta_l^{(j)}},\qquad j=0,\ldots,q,
  \la{11}
  \end{equation}
  where  $\beta_l^{(j)},\ l\geq1$ are i.i.d. random variables given
  by \refm[8],\refm[10]. Then
  \ber
  R(\overrightarrow{L^{\star}},\overrightarrow{N^{\star}})&=&c_N^{-1}\prod_{j=0}^{q}
  {{\frac{\left(S_j^{\star}\right)^{L_j^{\star}}}{{(L_j^{\star})}!}}
  \exp\left(\delta_jN_j^{\star}\right)\PP(Y^{(j)}=N_j^{\star})},
  \la{12}
  \ena
   where
   \be
   N_0^{\star}:=N-\sum_{j=1}^{q}{N_j^{\star}}.
   \la{406}
   \end{equation}
  \end{lemma}
  \nin {\bf Proof:}\quad  It follows from \refm[8]-\refm[11] and the multinomial enumeration formula
  for the number of ways to distribute $L_j^{\star}$ numbered balls over $L_j^{\star}$
   numbered urns, that
  \ber
  \prod_{j=0}^{q}{\PP(Y^{(j)}=N_j^{\star})}=\left(\prod_{j=0}^{q}{\frac{{(L_j^{\star})}!}{\left(S_j^{\star}\right)^{{L_j^{\star}}}}}
  {\exp\left({-\delta_jN_j^{\star}}\right)}\right)\sum_{\eta\in\om}\left(\ds{\prod_{k=1}^{N}{\frac{a_k^{n_k}}{{n_k}!}}}
  \mathbf{1}_{A}\right).
  \ena
  \nin By \refm[14],\refm[15] this gives \refm[12].
  $\blacksquare$\\
  Note that in contrast to  the standard one parameter scheme of Khintchine's method,
  the representation \refm[12] is based on $(q+1)$ free parameters
  $\delta_1,\ldots,\delta_q.$
   \begin{remark}
   One can see that \refm[12] is a  generalization of the representation formula for the total number of components which was
   obtained in \cite{EG}.
   \end{remark}
   \subsection{The choice of free parameters}
\ \ \    Firstly, we set in \refm[12]
$\delta_0=\delta_1=\ldots=\delta_q:=\delta$ and
 apply  the
 stratifications $$\overrightarrow{S}=(S_0,\ldots,S_q),\quad
 \overrightarrow{E}=(E_0,\ldots, E_q), \quad
 \overrightarrow{V}=(V_0,\ldots,V_q),$$ as defined in \refm[19],
  of the  quantities $S:=\sum_{k=1}^{N}{a_ke^{-\delta k}},$\ $E:=\sum_{k=1}^{N}{ka_ke^{-\delta k}}$
 and $V:=\sum_{k=1}^{N}{k^2a_ke^{-\delta k}}$ respectively.
 Next, we take in \refm[300] the vectors $\overrightarrow{L}=\big(L_0,\ldots,L_q\big)$ and $\overrightarrow{N}=\big(N_1,\ldots,N_q\big)$
  in the form
  \ber
    \non & & L_j=S_j+x_j\sqrt{S_0},\qquad j=0,\ldots,q,\\
   & &N_j=E_j+x_{j+q}\sqrt{V_0},\quad j=1,\ldots,q ,
  \la{20}
  \ena
   where $x_j, \  j=0,\ldots, 2q$\  are arbitrary reals. We also set $N_0=N-N_1$ and   adopt the
   following notation:
    \ber
    x_j^{\circ}=
  \begin{cases}
    x_j, & \text{if}\ j\in\{0,q\},\\
    x_{j}-x_{j+1},  & \text{if}\ 1\leq j\leq q-1,\\
    -x_{q+1}, & \text{if}\ j=q+1,\\
    x_{j-1}-x_{j}, & \text{if}\ q+2\leq j\leq 2q,\\
    x_{2q}, &\text{if}\ j=2q+1.
  \end{cases}
  \la{129}
  \ena
  It follows from \refm[129] that $\ds\sum_{j=q+1}^{2q+1}{x_j^{\circ}}=0$. By \refm[401] and \refm[20] we then have
  \ber
  L_j^{\star}=S_j^{\star}+x_j^{\circ}\sqrt{S_0}, & j=0,\ldots,q
  \la{119}
  \ena
  and
  \ber
  N_j^{\star}= E_j^{\star}+x_{j+q+1}^{\circ}\sqrt{V_0}, & & j=1,\ldots,q.
  \la{120}
  \ena
  Note that \refm[406], \refm[120] together with \refm[19] imply that
  \be
  N_0^{\star}=N-\sum_{j=1}^{q}{\left(E_j^{\star}+x_{q+j+1}^{\circ}\sqrt{V_0}\right)}=N-E_1+x_{q+1}^{\circ}\sqrt{V_0}.
  \la{93}
  \end{equation}
  We start by investigating the asymptotic behavior of
the probability
$R(\overrightarrow{L^{\star}},\overrightarrow{N^{\star}})$ given
by \refm[12],
   when  $N\rightarrow\infty$ and $\delta$ is fixed.  We adopt from \cite{frgr1}
    the following
  representation \refm[21] of the partition function $ c_N$ defined in \refm[80].
   Consider independent random variables $\xi_k,$\
  s.t.\ $k^{-1}\xi_k\sim Po(a_ke^{-\sigma k}),\ k\geq 1,$ where $ \sigma>0$ is a free parameter, and let
  $\Z=\ds\sum_{k=1}^{N}{\xi_k}.$ Then
   \be
    c_N=\PP(\Z=N)\exp{\left(N\sigma+S(\sigma)\right)},\quad N\geq 1,
  \la{21}
  \end{equation}
  where $S=S(\sigma)=\ds\sum_{k=1}^{N}{a_ke^{-\sigma k}}=\sum_{j=0}^q{S_j^{\star}(\sigma)}$. We set $\sigma=\delta$ .
  The forthcoming asymptotic analysis relies on the following fact that features the expansive case \refm[95] considered.
Let $\delta=\delta_N\to 0,$ as $N\to \infty$ and let the
stratification points be as in Theorem 4.1, i.e. $M_j=u_jr_N, \
j=0,\ldots,q, \ 0=u_0<u_1<\ldots<u_q$. Then
  \be
  \ds\lim_{N\rightarrow\infty}{S_j^{\star}}=\infty,\quad S_j^{\star}=O\left(S_0\right),\quad  N\rightarrow\infty,
  \quad j=0,\ldots,q.
  \la{sh}
  \end{equation}
 This is proved with the help of the integral test,
as it is detailed in
 \refm[38] below.
  We will see  from the proof of Theorem 4.6 that \refm[sh] fails for the asymptotics
  of   component counts
   of small sizes.

 \begin{lemma}
  Suppose \refm[sh] holds. Then for $L_0^{\star},\ldots,L_q^{\star}$ given by \refm[119] and for any \\
  $N_0^{\star},\ldots,N_q^{\star}\in\mathbb{N}:\ds\sum_{j=0}^q {N_j^{\star}}=N,$
  \ber
   R(\overrightarrow{L^{\star}},\overrightarrow{N^{\star}})&\sim & \frac{W}{{(2\pi)}^{\frac{q+1}{2}}}
  \left(\prod_{j=0}^{q}{\frac{1}{\sqrt{S_j^{\star}}}}\right)\exp\left({-\half
  \sum_{j=0}^{q}{{\left(x_j^{\circ}\right)}^2\frac{S_0}{S_j^{\star}}}}\right),\quad \  N\to \infty,
  \la{22}
  \ena
  where
  \be
  \ds W=\frac{\prod_{j=0}^{q}{\PP(Y^{(j)}=N_j^{\star})}}{\PP(Z_N=N)}.
  \la{121}
  \end{equation}
  \end{lemma}
  \nin {\bf Proof:} We have from \refm[12] and \refm[21]
  \ber
 \non R(\overrightarrow{L^{\star}},\overrightarrow{N^{\star}})&=&W\cdot\left(\prod_{j=0}^{q}\frac{\left(S_j^{\star}\right)^{{L_j^{\star}}}}
  {(L_j^{\star})!}\right)\exp{(\delta N)}\exp{\big(-(\delta N+S)\big)}\\&=& W\exp{\left(-S\right)}\prod_{j=0}^{q}
  {\frac{\left(S_j^{\star}\right)^{{L_j^{\star}}}}{(L_j^{\star})!}}.
  \la{24}
  \ena
  Next, we apply Stirling's asymptotic formula to estimate ${(L_j^{\star})!}$ under \refm[119] and the condition \refm[sh]:
  \ber
  \non({L_j^{\star}})!&\sim&\sqrt{2\pi L_j^{\star}}\;\left(L_j^{\star}\right)^{{L_j^{\star}}}\exp{(-L_j^{\star})}
  \sim\sqrt{2\pi S_j^{\star}}\:\left(S_j^{\star}\right)^{{L_j^{\star}}}\left(1+ \frac{ x_j^{\circ}\sqrt{S_0}}
  {S_j^{\star}}\right)^{L_j^{\star}}\exp{(-L_j^{\star})}\\ \non  &=& \sqrt{2\pi S_j^{\star}}\:
  \left(S_j^{\star}\right)^{{L_j^{\star}}}\exp{\left(L_j^{\star}\log\big(1+\frac{ x_j^{\circ}\sqrt{S_0}}
  {S_j^{\star}}\big)\right)}\exp(-L_j^{\star})\\
   &\sim &\sqrt{2\pi S_j^{\star}}\:\left(S_j^{\star}\right)^{{L_j^{\star}}}e^{-S_j^{\star}}\exp{\left(\half
   {\left(x_j^{\circ}\right)}^2\frac{S_0}{S_j^{\star}}\right)},\qquad
   N\to\infty, \quad j=0,\ldots,q.   \la{sh1}
  \la{25}
  \ena
  Substituting in \refm[24] the preceding asymptotic expansions yields \refm[22]
  $\blacksquare$\\

  We will follow  the principle of Khintchine's method  that the
free parameter
 $\delta$ should be coupled with $N$ so that to make
  the probabilities in the RHS of \refm[121] large. Namely, (see also  \cite{frgr1}, \cite{frgr2})  we choose $\delta=\delta_N$
  as the solution of the equation:
  \be
  \E Z_N=\sum_{k=1}^{N}{ka_ke^{-\delta k}}=N.
  \la{16}
  \end{equation}
   We see from \refm[19] that \refm[16] implies $E_0^{\star}=E_0=N-E_1$. So, we obtain from  \refm[93] and \refm[129] that
  \be
  N_0^{\star}=N-N_1=E_0^{\star}+x_{q+1}^{\circ}\sqrt{V_0},
  \la{160}
  \end{equation}
  which enables us to rewrite \refm[120] in a unified way:
\ber
  N_j^{\star}= E_j^{\star}+x_{j+q+1}^{\circ}\sqrt{V_0}, & j=0,\ldots,q.
  \la{920}
  \ena
  It is easy to see that if $a_k,\ k\geq 1$ are positive, the equation \refm[16] has a unique solution for any $N\geq1$.\\

  Our next goal will be the establishment of local limit theorems needed for the proof of
  Theorem 4.1..

 \subsection{The local limit theorems for $\yn,\ j=0,\ldots,q$.}
  \begin{lemma}\ (The local limit theorems for $\yn,\ j=0,\ldots,q$).\\
  Let \refm[95] hold and let $M_j=[u_jr_N],
  j=0,\ldots,q,$ where $ 0=u_0<u_1<\ldots<u_q, \    \ M_{q+1}=N+1$ and $[\bullet]$ is
  the integer part of a number. Then for
  $\overrightarrow{L^{\star}},\overrightarrow{N^{\star}}$ as in \refm[119],\refm[120]
   and for $\delta$ given  by
  \refm[16],
  \ber
    \PP(\yn=N_j^{\star})\sim\frac{1}{\sqrt{2\pi \Var{\yn}}}\exp{\left(-\frac{\kappa_j^2}
    {2}\right)},
   \quad N\rightarrow\infty,\quad j=0,\ldots,q,
  \la{29}
  \ena
  where $$\kappa_j=\lim_{N\to\infty}
  \frac{\E{\yn}-N_j^{\star}}{\sqrt{\Var{\yn}}}, \ \ \ j=0, \ldots,q$$ are positive constants
  calculated in \refm[50]  below.
  \end{lemma}
 \nin {\bf Proof:}
  We denote by $\varphi^{(j)}$ the characteristic functions of the random variable
  $\yn,$ to obtain
  \be
   \PP(\yn=N_j^{\star})={(2\pi)}^{-1}\int_{-\pi}^{\pi}{\wt e^{-iN_j^{\star}t}dt}:= {(2\pi)}^{-1}I^{(j)},
   \qquad j=0,\ldots, q.
  \la{30}
  \end{equation}
   We will focus now on the asymptotics, as $N\rightarrow\infty$, of the integrals $I^{(j)},\ j=0,\ldots,q$.
   For any $0<|t_{0,j}|<\pi$  the integral $\ii$ can be  written in the form:
   \be
   \ii=I_1^{(j)}+I_2^{(j)},\qquad j=0,\ldots q,
   \la{32}
   \end{equation}
   where $\ii_1=\ii_1(t_{0,j})$ and $\ii_2={\ii_2}(t_{0,j})$ are integrals of the
   integrand of $\ii,$ taken over
   the sets $[-t_{0,j},t_{0,j}]$ and $[-\pi,-t_{0,j}]\bigcup[t_{0,j},\pi]$,
   respectively.
   Using the technique of \cite{frgr1},\cite{frgr2}, we will show that for an
   appropriate choice of $0<t_{0,j}=t_0(N,j)\to 0,\ N\to \infty,\ j=0,\ldots,q$,
   the main contribution to $\ii$, as $N\rightarrow\infty$, comes from $\ii_1$, i.e. from a specially
constructed neighborhood of $t=0$. Firstly, observe that
   \be
   \wt={\left(\wz(t)\right)}^{L_j^{\star}},\quad t\in\mathbb{R},\quad
   j=0,\ldots,q,
   \la{33}
   \end{equation}
  where $\wz$ is the characteristic function of the random variable
  $\beta_j:=\beta_1^{(j)}$ defined by \refm[8],\refm[10].
  To derive the asymptotics
  of the integral $I_1^{(j)},$ under   the above choice of $t_{0,j},$ we will
   look for the approximation of $\wt$,\ as $t\rightarrow 0$ and $N\to \infty.$  For this purpose  we need the
  asymptotic expressions, as $N\to\infty$, for $\E{\beta_j},\;\E{\beta_j^2} $ and
  $\E{\beta_j^3}$. We have
  \ber
  \E{\beta_j}=\frac{E_j^{\star}}{S_j^{\star}},\qquad \E{\beta_j^2}=\frac{V_j^{\star}}{S_j^{\star}},\qquad
  \E{\beta_j^3}=\frac{H_j^{\star}}{S_j^{\star}},\qquad j=0,\ldots,q,
  \la{35}
  \ena
  where $H_j^{\star}=H_j^{\star}(\delta,M_j,M_{j+1})=\ds\sum_{k=M_j}^{M_{j+1}-1}{k^3a_ke^{-\delta k}},\ j=0,\ldots,q$.
  It follows from  the definitions of $S_j^{\star},E_j^{\star},V_j^{\star}$ and $H_j^{\star}$ that in the case considered
  the problem reduces to estimation of sums of\ the form
  \be
   \ds\sum_{k=M_j}^{M_{j+1}-1}{a_kk^re^{-\delta k}},\quad a_k\sim Ck^{p-1},\ k\rightarrow\infty,\quad
   p>0,\quad C>0,\quad  r\geq 0,
  \la{37}
  \end{equation}
  when $N\rightarrow\infty$,   $\delta$ is given as the  solution of \refm[16] and
  $M_j, j=0,\ldots,q+1$  are as
  in the statement of the Lemma. The asymptotic solution of \refm[16] was obtained in \cite{frgr2}:
  \be
   \delta=\delta_N\sim hN^{-\frac{1}{p+1}}=(r_N)^{-1},\quad N\rightarrow\infty,
  \la{36}
  \end{equation}
  where $h, r_N$ are as in \refm[h]. Thus,  $ M_j\delta\to u_j, \ N\to \infty,\ j=0,\ldots,q.$
  In view of (\ref{95}), it is convenient to write $a_k=Ck^{p-1}G(k)$, where the function $G$ is s.t. $\ds\lim_{k\rightarrow\infty}{G(k)}=1$.
   Then, applying the integral test (=Euler summation formula), we get from \refm[36]
  \ber
   \sum_{k=M_j}^{M_{j+1}-1}{a_kk^re^{-\delta k}}& \sim  &C\int_{M_j}^{M_{j+1}-1}{x^{p+r-1}G(x)e^{-\delta\non x}dx}
   =C\delta^{-p-r}\int_{M_j\delta}^{(M_{j+1}-1)\delta}{x^{p+r-1}G{\left(\frac{x}{\delta}\right)}e^{-x}dx}\\ \non &\sim&
   \ds C\delta^{-p-r}\int_{u_j}^{u_{j+1}}{x^{p+r-1}G{\left(\frac{x}{\delta}\right)}e^{-x}dx}\sim\delta^{-p-r}
   C\int_{u_j}^{u_{j+1}}{x^{p+r-1}e^{-x}dx}\\ &:=&
   \delta^{-p-r} f_j(p+r-1),\quad N\rightarrow\infty,\quad j=0,\ldots,q,
   \quad r\geq0,\ p>0,
  \la{38}
  \ena
  where we set $u_{q+1}=\infty$. Note that by our notation \refm[26],
  \ber
 f_j(s)=b_s(u_j)-b_s(u_{j+1}),\quad s>-1,\quad j=0,\dots,q .
  \la{94}
  \ena
  We also obtain from \refm[38] and \refm[35],
  \be
  \E{\beta_j^r}=O(\delta^{-r}),\quad N\rightarrow\infty,\quad r>0,\quad
  j=0,\ldots,q
  \la{39}
  \end{equation}
 and observe that
  \be
  f_j(p+1) f_j(p-1) -(f_j(p))^2>0, \quad f_j(p+2)f_j(p-1)-f_j(p+1)f_j(p)>0,\quad j=0,\dots,q.
 \la{cd}
 \end{equation}
 The first of these inequalities follows immediately from the Cauchy-Schwarz inequality, while  the second can
 be derived by substituting $f_j(p-1)$ from the first one and then applying again the Cauchy-Schwarz inequality.\\
  To arrive at  the required asymptotic formula for $\varphi^{(j)},$  we first write for a fixed $N$,
  \be
   \wz(t)=1+i\E{\beta_j}t-\half \E{\beta_j^2}t^2+O\big(\E{\beta_j^3}t^3\big),\quad t\rightarrow0,\quad j=0,\ldots,q,
  \la{41}
  \end{equation}
  and then  couple  $t$ with $N$ by setting
   \be
   t_{0,j}=\delta {(L_j^{\star})}^{-\half}\log^2{\delta}=
   O( \delta^{1+\frac{p}{2}}\log^2\delta),\quad
   j=0,\ldots,q,\quad N\to \infty,
   \la{47}
   \end{equation}
   where the last equality follows from \refm[38] and \refm[119]. Consequently,
   \be
   \left(\E{\beta_j^r}\right)t_{0,j}^r=O(\delta^\frac{rp}{2}\log^{2r}\delta)\to 0,
   \quad N\rightarrow\infty, \quad j=0,\ldots,q,\quad r> 0.
   \la{327}
   \end{equation}
    Now \refm[33] and  \refm[41] together with \refm[119], \refm[120] give
   \ber
   \non & &\wt e^{-iN_j^{\star}t}={\left(1+i\E{\beta_j}t-\half \E{\beta_j^2}t^2+ O\big(\E{\beta_j^3}t^3\big)\right)}
   ^{L_j^{\star}}e^{-iN_j^{\star}t}\\ \non &\sim & \exp{\left(it\left(L_j^{\star}\E{\beta_j}-N_j^{\star}\right)
   -\half L_j^{\star}\E{\beta_j^2}t^2+\half L_j^{\star}t^2{\big(\E{\beta_j}\big)}^2+L_j^{\star}O\big(\E{\beta_j^3}t^3\big)
   \right)}\\ \non &= & \exp{\left(it\left(\E{\yn}-N_j^{\star}\right)-\half t^2 \Var{\yn}+L_j^{\star}
   O\big(\E{\beta_j^3}t^3\big)\right)},\\  & &  |t|\leq t_{0,j},\quad N\rightarrow\infty,\quad j=0,\ldots,q.
   \la{42}
   \ena
     In view of \refm[39] and  \refm[cd], the  choice \refm[47] of $t_{0,j}$ provides
   \be
     \lim_{N\rightarrow\infty}t_{0,j}^2{\Var{\yn}}=\infty,\quad
     j=0,\ldots,q
   \la{43}
   \end{equation}
   and
    \be
    \lim_{N\rightarrow\infty}{L_j^{\star}t_{0,j}^3\E{\beta_j^3}}=0,\quad j=0,\ldots,q.
   \la{44}
   \end{equation}
   \begin{remark}
   By \refm[39] and \refm[cd],
   \be
  \Var{\yn}=L_j^{\star}\left(\E{\beta_j^2}-{\left(\E{\beta_j}\right)}^2\right)= O\big(\delta^{-2}L_j^{\star}\big),
  \quad N\rightarrow\infty,\quad j=0,\ldots,q
   \la{45}
   \end{equation}
   and
   \ber
   \la{73}
   \non\E{\left(\yn-\E{\yn}\right)}^3&=&L_j^{\star}\E{\left(\beta_j-\E{\beta_j}\right)}^3=
   L_j^{\star}\left(\E{\beta_j^3}-3\E{\beta_j^2}\E{\beta_j}+2{\big(\E{\beta_j}\big)}^3\right)=
   O\big(\delta^{-3}L_j^{\star}\big),\\ & & N\rightarrow\infty,\quad
   j=0,\ldots,q.
   \ena
   Hence, the following weaker (=the third moment is not absolute) form of
   Lyapunov's sufficient condition for the convergence to the normal law
   holds for the
   sums $\yn,\ j=0,\ldots,q$:
   \be
   \lim_{N\rightarrow\infty}{\frac{\E{\left(\yn-\E{\yn}\right)}^3}{{\left(\Var{\yn}\right)}^{\frac{3}{2}}}}=0
   ,\quad  j=0,\ldots,q.
   \la{74}
   \end{equation}
   This explains the existence of $t_{0,j}$ that provides \refm[43] and \refm[44].
  The phenomenon described above is typical in applications of Khintchine's method (see \cite{fr},\cite{frgr2},\cite{GS}).
   \end{remark}
    To continue the asymptotic expansion \refm[42], we obtain from \refm[38], \refm[35] and \refm[cd]
    \ber
   \non E_j^{\star}\sim f_j(p)\delta^{-p-1},\quad V_j^{\star} &\sim& f_j(p+1){\delta}^{-p-2},\quad \Var{\yn}\sim
    w_j\delta^{-p-2},\quad S_j^{\star}\sim f_j(p-1){\delta}^{-p},\\ & &\quad N\rightarrow\infty,\quad j=0,\ldots,q,
    \la{75}
    \ena
    where
   \be
   w_j=f_j(p+1)-\frac{f_j^2(p)}{f_j(p-1)}>0,\quad j=0,\ldots,q.
   \la{51}
   \end{equation}
 \nin The  asymptotic relations below are the consequence of \refm[75], \refm[119] and
 \refm[120]:
   \ber
   \non& &\frac{\E{\yn}-N_j^{\star}}{\sqrt{\Var{\yn}}}=\frac{\frac{E_j^{\star}}{S_j^{\star}}L_j^{\star}-N_j^{\star}}
   {\sqrt{\Var{\yn}}}=x_j^{\circ}\frac{E_j^{\star}\sqrt{S_0}}{S_j^{\star}\sqrt{\Var{\yn}}}
   -x_{j+q+1}^{\circ}\frac{\sqrt{V_0}}{\sqrt{\Var{\yn}}}\\ \non &\sim& x_j^{\circ}\frac{f_j(p)}{f_j(p-1)}
   \sqrt{\frac{{f_0(p-1)}}{w_j}}-x_{j+q+1}^{\circ}\sqrt{\frac{f_0(p+1)}{w_j}}:
   =\kappa_j(x_j^{\circ},x_{j+q+1}^{\circ})
   =\kappa_j,\\ & & N\rightarrow\infty, \quad\  j=0,\ldots,q.
   \la{50}
   \ena
   Now we are in a position to evaluate the integrals $\ii_1,\ j=0,\ldots,q $. By virtue of \refm[42]-\refm[44] and
   \refm[50],
   \ber
   \non \ii_1 &=& \int_{-t_{0,j}}^{t_{0,j}}{\wt e^{-iN_j^{\star}t}dt}\sim\int_{-t_{0,j}}^{t_{0,j}}
   {\exp{\left(it\big(\E{\yn}-N_j^{\star}\big)-\half t^2\Var{\yn}\right)}dt}\\ \non &=& \exp\left({-\frac{{(\E{\yn}-
   N_j^{\star})}^2}{2\Var{\yn}}}\right)\int_{-t_{0,j}}^{t_{0,j}}{\exp\left[-\half{\left(t\sqrt{\Var{\yn}}-i
   \frac{\E{\yn}-N_j^{\star}}{\sqrt{\Var{\yn}}}\right)}^2\right]dt}\\\non &\sim &\frac{1}{\sqrt{\Var{\yn}}}\exp
   \left({-\frac{{(\E{\yn}-N_j^{\star})}^2}{2\Var{\yn}}}\right)\int_{-\sqrt{\Var{\yn}}t_{0,j}}^{\sqrt{\Var{\yn}}t_{0,j}}
   {\exp\left[-\half{\left(t-i\kappa_j\right)^2}\right]dt}\\ \non &\sim &\frac{1}{\sqrt{\Var{\yn}}}\exp
   \left({-\frac{{(\E{\yn}-N_j^{\star})}^2}{2\Var{\yn}}}\right)\int_{-\infty}^{\infty}{\exp\left[-\half{\left(t-i\kappa_j
   \right)^2}\right]dt}\\ &\sim & \frac{\sqrt{2\pi}}{\sqrt{\Var{\yn}}}\exp{\left({-\frac{\kappa_j^2}{2}}\right)},
   \quad N\rightarrow\infty, \quad j=0,\ldots,q.
   \la{48}
   \ena
    Next we turn  to the estimation, as $N\rightarrow\infty$, of the integrals $\ii_2,\ j=0,\ldots,q$. We start with
    \ber
     \non & & |\ii_2|:=2\left|\int_{t_{0,j}}^{\pi}{\wt e^{-iN_j^{\star}t}dt}\right|\leq2\int_{t_{0,j}}^{\pi}{|\wt|dt}
     =4\pi\int_{t_{1,j}}^{\half}{\vert\varphi^{(j)}(2\pi t)\vert dt}\\ &= & 4\pi\int_{t_{1,j}}^{\half}
     {\left|\sum_{k=M_j}^{M_{j+1}-1}{a_ke^{(2\pi it-\delta)k}}\right|^{L_j^{\star}}{\left({S_j}^{\star}\right)}
     ^{-L_j^{\star}}dt},\quad j=0,\ldots,q,
    \la{83}
     \ena
     where we set $t_{1,j}=\frac{t_{0,j}}{2\pi}>0,\ j=0,\ldots,q$. Denote
   \be
    g^{(j)}(t)=\vert\varphi_1^{(j)}(2\pi t)\vert, \ t\in [t_{1,j},1/2],
    \quad j=0,\ldots,q.
    \la{148}
   \end{equation}
   It is easy to see from the definition of the random variables $\beta_j$ that
  $g^{(j)}(1/2)<1,\ j=0,\ldots,q.$ Since the  $\beta_j,\ j=0,\ldots,q$
   have lattice distributions with
  span $1,$ this implies (see \cite{dur}, p.131, \cite{shir}, p.286)
  that for a fixed $N,$ we have
  $g^{(j)}(t)<1, \ t\in [t_{1,j},1/2],\ j=0,\ldots,q.$ Moreover,
   $g^{(j)}(t_{1,j})\to 1, \quad N\to \infty,$
  because  $t_{1,j}\to 0, \quad N\to \infty.$ We then  conclude that
  for sufficiently large $N $ and any $0<t\le 1/2,$
  \be
   g^{(j)}( t_{1,j})\ge g^{(j)}( t),\quad j=0,\ldots,q.
   \la{hv}
 \end{equation}
 By the same argument,
 \ber
  \non\varphi_1^{(j)}(2\pi t_{1,j})-1&=&{\left(S_j^{\star}\right)}^{-1}\sum_{k=M_j}^{M_{j+1}-1}{a_ke^{-\delta
     k}(e^{2\pi ikt_{1,j}}-1)}\to 0,\\ & & \ N\to \infty,\quad
     j=0,\ldots,q.
 \la{147}
  \ena

\nin  This and \refm[20] enable us  to write
   \ber
   \non {\left(g^{(j)}(t_{1,j})\right)}^{L_j^{\star}} &=&\left\vert 1+{\left(S_j^{\star}\right)}^{-1}
   \sum_{k=M_j}^{M_{j+1}-1}{a_ke^{-\delta k}(e^{2\pi ikt_{1,j}}-1)}\right\vert^{L_j^{\star}}\\ \non
   &\sim& \exp{\left(-2\big(1+O(\delta^{\frac{p}{2}})\big)\sum_{k=M_j}^{M_{j+1}-1}{a_ke^{-\delta k}
   \sin^2{(\pi kt_{1,j})}}\right)},\\ & & N\rightarrow\infty,\quad j=0,\ldots,q.
   \la{55}
     \ena
     Denote
     \be
      D_j=2\sum_{k=M_j}^{M_{j+1}-1}{a_ke^{-\delta k}\sin^2{(\pi kt_{1,j})}}
     \quad  j=0,\ldots,q
     \la{56}
     \end{equation}
   and let $\tilde{u}_j\in(u_j,u_{j+1}),\ j=0,\ldots,q$. Then we have
     \ber
    D_j\non &\geq &2\sum_{k=u_jN^{\frac{1}{p+1}}}^{\tilde{u}_jN^{\frac{1}{p+1}}}{a_ke^{-\delta k}\sin^2{(\pi kt_{1,j})}}
     \geq2\left(\min_{k\in[u_jN^{\frac{1}{p+1}},\tilde{u}_jN^{\frac{1}{p+1}}]}({a_ke^{-\delta
     k}})\right)\sum_{k=u_jN^{\frac{1}{p+1}}}^{\tilde{u}_jN^{\frac{1}{p+1}}}{\sin^2{(\pi kt_{1,j})}}=\\ &&
     O\left(\delta^{-p+1}\right)\sum_{k=u_jN^{\frac{1}{p+1}}} ^{\tilde{u}_jN^{\frac{1}{p+1}}}{\sin^2{(\pi kt_{1,j})}},\quad
     N\rightarrow\infty,\quad j=0,\ldots,q.
     \la{58}
     \ena
     Denote $\epsilon_j=\tilde{u}_j-u_j>0,\ j=0,\ldots,q$. To estimate the last sum in
     \refm[58],
     we employ the following inequality  from \cite{fr}:
     \be
     2\sum_{k=l}^{l+m-1}{\sin^2{(\pi tk)}}\geq\frac{m}{2}\min\{1,(tm)^2\},\quad |t|\leq\half,\quad m\geq2,\quad l\geq1.
     \la{57}
      \end{equation}
  We apply \refm[57] with $l=u_jN^{\frac{1}{p+1}}$ and $m=\epsilon_j N^{\frac{1}{p+1}}+1$ to get from \refm[47] and
  \refm[36],
   \ber
    \non & &2\sum_{k=u_jN^{\frac{1}{p+1}}}^{\tilde{u}_jN^{\frac{1}{p+1}}}{\sin^2(\pi kt_{1,j})}
    \geq\half\epsilon_jN^{\frac{1}{p+1}}\min\{1,(\epsilon_jN^{\frac{1}{p+1}}t_{1,j})^2\}\\
    =& &O\left(\delta^{p-1}\right)\log^4{\delta},\quad  N\rightarrow\infty,\quad j=0,\ldots,q.
    \la{59}
    \ena
  Now we deduce from \refm[58], \refm[59], \refm[83] and \refm[hv] that
  \be
  \vert I_2^{(j)}\vert\le O(\exp{(-\log^4{\delta}})), \quad N\to\infty.
  \la{vk}
   \end{equation}
 Comparing this with \refm[48] gives
    \be
     \ii_2=o(\ii_1),\quad N\rightarrow\infty,\quad j=0,\ldots q,
    \la{60}
    \end{equation}
    which together with \refm[50], \refm[48] proves the lemma.\quad $\blacksquare$
 \subsection{Completion of the proofs of Theorems 4.1,\ 4.3 and 4.6}

\ \ \     To complete the asymptotic analysis of the probability
    $R(\overrightarrow{L^{\star}},\overrightarrow{N^{\star}})$ in the vicinity of $r_N$, it
    remains to derive the asymptotics of $\PP(Z_N=N),$ as $N\to \infty.$  It was found  in \cite{frgr1} that
    \be
   \PP(Z_N=N)\sim\frac{1}{\sqrt{2\pi{\Var{Z_N}}}},\quad N\rightarrow\infty,
   \la{61}
   \end{equation}
 \nin where
  \be
   \Var{Z_N}=\sum_{k=1}^{N}{k^2a_ke^{-\delta k}}\sim\Gamma(p+2)\delta^{-(p+2)},\quad N\rightarrow\infty.
   \la{62}
  \end{equation}
 \nin Substituting  in \refm[22] the asymptotic expressions \refm[29],\refm[75],\refm[61] and
 \refm[62] we  obtain
  \ber
  \non & &R(\overrightarrow{L^{\star}},\overrightarrow{N^{\star}})\sim\frac{\sqrt{\Var{Z_N}}}{{(2\pi)}^{q+\half}}
  \left(\prod_{j=0}^{q}\ds\frac{1}{\sqrt{S_j^{\star}\Var{Y^{(j)}}}}\right)\exp{\left[-\half\left(\sum_{j=0}^{q}
  {{(x_j^{\circ})}^2\frac{f_0(p-1)}{f_j(p-1)}}+\kappa_j^2\right)\right]}\\ \non&\sim& \delta^{(q(p+1)+p/2)}
   \frac{\sqrt{\Gamma(p+2)}}{{(2\pi)}^{q+\half}}\left(\prod_{j=0}^{q}{\frac{1}{\sqrt{w_jf_j(p-1)}}}\right)
  \exp{\left[-\half\left(\sum_{j=0}^{q}{{(x_j^{\circ})}^2\frac{f_0(p-1)}{f_j(p-1)}}+\kappa_j^2\right)\right]},\\  & &
  \quad N\rightarrow\infty.
  \la{065}
  \ena
 \nin Next we employ the scaled quantities $\hat{\nu}_j^{\star}, \ \hat{K}^{\star}_j$
  constructed in \refm[68]
  and denote
 \ber
 \non
 \overrightarrow{\hat{\nu}^{\star}}=(\hat{\nu}_0^{\star},\ldots,\hat{\nu}_q^{\star}),\qquad\overrightarrow{\hat{K}^{\star}}
 =(\hat{K}_1^{\star},\ldots,\hat{K}_q^{\star}),\\ \overrightarrow{x^\circ}=(x_0^\circ,\ldots,x_q^\circ),\qquad
  \overrightarrow{y^\circ}=(x_{q+2}^\circ,\ldots,x_{2q+1}^\circ).
 \la{169}
 \ena
 Note that  $\hat{K_0^{\star}}=-\sum_{j=1}^{q}{\hat{K}_j^{\star}},$ by \refm[16]. We get from \refm[065]
 \ber
 \non &&\PP \left(\overrightarrow{\hat{\nu}^{\star}}=\overrightarrow{x^{\circ}},
 \overrightarrow{\hat{K}^{\star}}=\overrightarrow{y}^{\circ}\right)\sim
 \\\non & &\frac{\sqrt{\Gamma(p+2)}}{{(2\pi)}^{\frac{2q+1}{2}}}{\left(\frac{f_0(p-1)}{S_0}
  \right)}^{\half{(q+1)}}{\left(\frac{f_0(p+1)}{V_0}\right)}^{\half{q}}
  \left(\prod_{j=0}^{q}{\frac{1}{\sqrt{w_jf_j(p-1)}}}\right)\times\\ & &\exp{\left[-\half\sum_{j=0}^{q}
  {\frac{{(x_j^{\circ})}^2}{f_j(p-1)}}\right]}\exp{\Bigg[-\half\sum_{j=0}^{q}\bigg(x_j^{\circ}{\frac{f_j(p)}{f_j(p-1)
 \sqrt{w_j}}-x_{q+j+1}^{\circ}
  \frac{1}{\sqrt{w_j}}}\bigg)^2\Bigg]},\quad N\rightarrow\infty,
 \la{200}
 \ena
 where $x_{q+1}^{\circ}=-\sum_{j=1}^{q}{x_{q+j+1}^{\circ}}$. \\
 \nin Each term of the sum in the second exponent depends on the two variables  $x_j^{\circ}$
 and $x_{q+j+1}^{\circ}$ only. Hence, under a fixed vector $\overrightarrow{y^{\circ}},$ the exponents in the RHS
 of \refm[200] factorize into a product of $q+1$ terms, each one depending on only one  of the $x_j^{\circ},
 \ j=0,\ldots, q.$ This proves the claim (i) of  Theorem 4.1.

 For the proof the claim  {\rm(ii)}  of  Theorem 4.1, we firstly need to deduce from \refm[200] the central limit theorem for
 the  $(2q+1)-$ dimensional random vector $(\overrightarrow{\hat{\nu}^{\star}},\overrightarrow{\hat{K}^{\star}})$.
 We see from the definition
 \refm[50] of $\kappa_j,\ j=0,\ldots,q$ that the expression in the product of the exponents in \refm[200] is a
 quadratic form of the coordinates of the  vector $(\overrightarrow{x^{\circ}},\overrightarrow{y^{\circ}})$. For a given $q\geq1$, we denote
  by $\Theta^{\star}(q)$ the $(2q+1)\times(2q+1)$ matrix of this quadratic form.
  Let now $c_j<d_j,\  j=0,\ldots,2q$ and
  \be
  (\overrightarrow{c},\overrightarrow{d})=
  (c_0,d_0)\times(c_1,d_1)\times\ldots\times(c_{2q},d_{2q}).
  \la{201}
  \end{equation}
  We also define the sets of discrete points
  \ber
  G_N^{(j)}=
  \begin{cases}
  \{z\in(c_j,d_j):S_j^{\star}+z\sqrt{\frac{S_0}{f_0(p-1)}}\in\mathbb{N}\},           & \text{if}\ \ 0\leq j\leq q,\\
  \{z\in(c_j,d_j):E_j^{\star}+z\sqrt{\frac{V_0}{f_0(p+1)}}\in\mathbb{N}\},           & \text{if}\ \ q+1\le j\leq 2q
  \end{cases}
  \la{205}
  \ena
  and $\overrightarrow{G_N}=G_N^{(0)}\times\ldots\times G_N^{(2q)}$.
  \begin{lemma}( The central limit theorem for the vector ($\overrightarrow{\hat{\nu}^{\star}},\overrightarrow
  {\hat{K}^{\star}}$)).\\
   As $N\rightarrow\infty,$ and $q\geq1$ is fixed,
  \ber
  \PP\left((\overrightarrow{\hat{\nu}^{\star}},\overrightarrow{\hat{K}^{\star}})\in(\overrightarrow{c},\overrightarrow{d})
  \right)\sim{(2\pi)}^{-\frac{2q+1}{2}}\sqrt{T_q}\int_{(\overrightarrow{c},\overrightarrow{d})}{\exp{\left(-\half
  \overrightarrow{z}^{T}\Theta^{\star}(q)\overrightarrow{z}\right)}d\overrightarrow{z}},
  \la{202}
  \ena
  where the matrix $\Theta^{\star}(q)$ is defined as above, and
  \be
  T_q=\Gamma(p+2)\left(\prod_{j=0}^{q}{\frac{1}{w_jf_j(p-1)}}\right).
  \la{203}
  \end{equation}
  \end{lemma}
  {\bf Proof:} We follow the known technique for passing from a local limit theorem  to an  integral one, that is exposed
  in detail in \cite{shir}, p.60 (see also \cite{dur}, p.80 and \cite{EG}). Summing \refm[200] over $(\overrightarrow{x^{\circ}},\overrightarrow{y^{\circ}})
  \in\overrightarrow{G_N}$, we obtain
  \ber
 \non & & \PP\left((\overrightarrow{\hat{\nu}^{\star}},\overrightarrow{\hat{K}^{\star}})\in(\overrightarrow{c},\overrightarrow{d})\right)=
  {(2\pi)}^{-\frac{2q+1}{2}}\sqrt{T_q}{\left(\frac{f_0(p-1)}{S_0}\right)}^{\half{(q+1)}}
 {\left(\frac{f_0(p+1)}{V_0}\right)}^{\half{q}}\\ & &\sum_{\overrightarrow{z}\in
   \overrightarrow{G_N}}{\left[\exp{\left(-\half\overrightarrow{z}^{T}\Theta^{\star}(q)\overrightarrow{z}\right)}
   \left(1+\epsilon(\overrightarrow{z},N)\right)\right]},\quad N\to\infty.
  \la{204}
  \ena
 \nin From the preceding asymptotic formulae  we derive the crucial fact  that in \refm[204],
 \be
 \sup_{\overrightarrow{z}\in(\overrightarrow{c},\overrightarrow{d})}{\epsilon(\overrightarrow{z},N)}\to 0,\quad N\to
 \infty.
 \la{1zx}
 \end{equation}

 \nin  Based on this property of uniform convergence, we treat the RHS of \refm[204] as a Riemann sum with the
 asymptotically equidistant spacings $\left(\vert G_N^{(j)}\vert\right)^{-1}\sim
  \sqrt{\frac{f_0(p-1)}{S_0}}\sim \delta^{\frac{p}{2}},\ j=0,\ldots,q$
  and
  $\left(|G_N^{(j)}|\right)^{-1}\sim \sqrt{\frac{f_0(p+1)}{V_0}}\sim \delta^{\frac{p+2}{2}},\ j=q+1,\ldots,2q,$
  as $N\to \infty.$ \qquad $\blacksquare$\\
  Next, we have to show that the main term in the RHS of \refm[204] is indeed the gaussian distribution. This is equivalent
  to proving that
  \be
  \det{\Theta^{\star}(q)}=T_q,\quad q\geq1,
  \la{305}
  \end{equation}
  where $T_q$ is as in \refm[203]. By the definition \refm[51] of $w_j,$
  \be
  \frac{1}{f_j(p-1)}+\frac{f_j^2(p)}{f_j^2(p-1)w_j}=\frac{f_j(p+1)}{f_j(p-1)w_j},\quad
  j=0,\ldots,q.
  \la{310}
  \end{equation}
  We first prove \refm[305] for $q=1$. Using \refm[310] we have from \refm[200]
  \ber
  \Theta^{\star}(1)=\left(  \begin{array}{cccc}
    \frac{f_0(p+1)}{f_0(p-1)w_0} &              0               & \frac{f_0(p)}{f_0(p-1)w_0}\\
    0                            & \frac{f_1(p+1)}{f_1(p-1)w_1} & -\frac{f_1(p)}{f_1(p-1)w_1}\\
    \frac{f_0(p)}{f_0(p-1)w_0}   & -\frac{f_1(p)}{f_1(p-1)w_1}  &\frac{1}{w_0}+\frac{1}{w_1}
     \end{array} \right).
  \la{309}
  \ena
  Now the claim  for $q=1$ follows
  from the identity $f_0(p+1)+f_1(p+1)=C\Gamma(p+2), \ p>0$ and some algebra.
  This shows that the random vector $(\hat{\nu}_0^{\star},\hat{\nu}_1^{\star},\hat{K}_1^{\star})$ weakly converges to the gaussian
  random vector with zero mean and the covariance matrix  $\left({\Theta^{\star}(1)}\right)^{-1}$.  This implies
  (see \cite{dur}, p.89, Theorem 2.6) that for $q=1$ the sequence of the distribution functions of the random vector
  $(\hat{\nu}_0^{\star},\hat{\nu}_1^{\star},\hat{K}_1^{\star})$ is tight under each
  $u_1>0$. Consequently, by the definition of tightness and the form of the stratifications considered, we deduce
  the tightness of the distribution functions of the random vector  $(\overrightarrow{\hat{\nu}^{\star}},
  \overrightarrow{\hat{K}^{\star}})$, for all $q>1$.
  \nin Finally, by Prohorov's theorem,  the limiting distribution of the above vector is the one of a  probability
   measure. This proves \refm[305] for $q\geq 1$.
  \vskip 0.2cm
  Now we derive the explicit form of the covariance matrix in Theorem 4.1.
  For a given $q\ge 1,$ let ${\left(\Theta^{\star}(q)\right)}^{-1}={\{\vartheta_{mk}^{\star}(q)\}}_{m,k=0}^{2q}
  $ be the covariance matrix of the limiting gaussian distribution in \refm[202].
   \begin{lemma}
   \begin{small}
   \ber
   \vartheta_{mk}^{\star}(q)=
  \begin{cases}
   f_m(p-1)\mathbf{1}_{m=k}-\frac{f_m(p)f_k(p)}{\Gamma(p+2)},             & \text{if}\ \ 0\leq m,k\leq q,\\
   f_{m}(p)\mathbf{1}_{m=k-q}-\frac{f_{m}(p)f_{k-q}(p+1)}{\Gamma(p+2)}, & \text{if}\ \ 0\leq m\leq q,\; q+1\leq k\leq 2q,\\
   f_{k}(p)\mathbf{1}_{k=m-q}-\frac{f_{k}(p)f_{m-q}(p+1)}{\Gamma(p+2)}, & \text{if}\ \ 0\leq k\leq q,\; q+1\leq m\leq 2q,\\
   f_{m-q}(p+1)\mathbf{1}_{m=k}-\frac{f_{m-q}(p+1)f_{k-q}(p+1)}{\Gamma(p+2)},       & \text{if}\ \ q+1\leq m,k\leq 2q.
  \end{cases}
   \la{covm}
   \ena
   \end{small}
   \end{lemma}
  {\bf Proof:} To understand the technique behind the inversion of the matrix $\Theta^{\star}(q)$, it is convenient to
  verify \refm[covm] first for $q=2.$ This is easy to do with the help  of  the identity $\sum_{j=0}^{q}{f_j(p+1)}
  =C\Gamma(p+2)$, which holds for any $q\geq1$. Taking into account the structure of the matrix $\Theta^{\star}(q),$
   the verification for  $q>2$ can be  done in the same way.\qquad
   $\blacksquare$\\
\nin This completes the proof of the statement \rm(ii) of Theorem
4.1..We note that by \refm[covm], all
 covariances    $\vartheta_{mk}^{\star}(q)$ between the scaled component counts $\hat{\nu}_j^{\star}$,
 as well as between masses $\hat{K}_j^{\star}$  are negative.

   The convergence of moments, as stated in {\rm(iii)} of Theorem 4.1 results from  \refm[1zx]. For the proof, it
 is needed only  to replace probabilities with moments in the argument given
 in \cite{shir}, p.61 (see also \cite{kal}, p.
 67).
   We see from \refm[covm] that for the proof of the claim ${\rm(iv)}$ of Theorem 4.1 one has to examine the behavior
 in $k\to \infty$ of the integrals $f_k(p),\ f_{k-q}(p+1)$ in the case
 when $u_j, \ j=0,\ldots, q$ are equidistant points and $q\to \infty.$
Now  Theorem 4.1 is completely proved.

   We  proceed  to the proof \ of  Theorem 4.3.
  First, we recall the scaled quantities  $\hat{\nu}_j,\ \hat{K}_j$ in \refm[104]
 and denote
 \ber
 \overrightarrow{\hat{\nu}}=(\hat{\nu}_0,\ldots,\hat{\nu}_q),\qquad \overrightarrow{\hat{K}}=(\hat{K}_1,\ldots,
 \hat{K}_q).
 \la{168}
 \ena
 It follows from the definition of $\nu_j$ and $K_j,\
 j=0,\ldots,q$ that
 \ber
 \nu_0=\nu_0^{\star},\quad\nu_j=\sum_{i=j}^{q}{\nu_i^{\star}},\quad K_j=\sum_{i=j}^{q}
 {K_i^{\star}},\quad j=1,\ldots,q.
 \la{222}
 \ena
  By Theorem 4.1, \refm[222] and \refm[covm], the  random vector $(\overrightarrow{\hat{\nu}}, \overrightarrow{\hat{K}})$
  converges weakly to the gaussian vector with zero mean and the covariance matrix $\Upsilon(q)={\{e_{mk}\}}_{m,k=0}^{2q}$
 given by:
 \begin{small}
  \ber
  e_{mk}=
  \begin{cases}
   f_{0}(p-1)-\frac{f_{0}^2(p-1)}{\Gamma(p+2)},                                 & \text{if}\  \ m=k=0,\\
   -\frac{f_{0}(p)b_{p}(u_k)}{\Gamma(p+2)},                                     & \text{if}\  \ m=0,\ 0<k\leq q,\\
   -\frac{f_{0}(p)b_{p+1}(u_{k-q})}{\Gamma(p+2)},                               & \text{if}\  \ m=0,\ q+1\leq k\leq 2q,\\
   b_{p-1}(u_{\max(m,k)})-\frac{b_{p}(u_m)b_{p}(u_k)}{\Gamma(p+2)},              & \text{if}\  \ 0<m,k\leq q,\\
   b_{p}(u_{\max(m, k-q)})-\frac{b_{p}(u_m)b_{p+1}(u_{k-q})}{\Gamma(p+2)},          & \text{if}\  \ 0<m\leq q,\; q+1\leq k\leq 2q,\\
   b_{p+1}(u_{\max(m,k)-q})-\frac{b_{p+1}(u_{m-q})b_{p+1}(u_{k-q})}{\Gamma(p+2)},& \text{if}\  \ q+1\leq m,k\leq 2q.
  \end{cases}
  \la{306}
  \ena
  \end{small}
  Applying now \refm[306] with $ 0<m,k\leq q$ gives the claim (i) of Theorem 4.3.

 \nin  For the proof of the claim {\rm(ii)} of the theorem,  we first recall that in the notation \refm[85], $\nu_j=\nu(r_Nu_j),
 \quad j=1,\ldots,q,$ with $r_N$ as in \refm[h]. Thus, letting $q=1$ and $u_1=u,$
 we rewrite the first part of \refm[104] as
 \be
 \hat\nu_1=\sqrt{f_0(p-1)}\ \ \frac{\nu(r_Nu)-S_1}{\sqrt{S_0}}
 \la{tos}
 \end{equation}
 to get with the help of \refm[38], \refm[94], \refm[h] and the notation \refm[86],
the desired relation \refm[relat] between $\hat{\nu}_j$ and the
fluctuations $\Delta(u_j) $ for $j=1,\ldots, q$.
 By virtue of the convergence stated in {\rm(i)} of Theorem 4.3, the fact that
 $Nr_N^{-1}\to \infty$, as $ N\to \infty$ and the continuity  of $Var \hat\nu_1 $ in
   $u,$
 we derive that $\hat\nu_1$ is stochastically equicontinuous (see e.g. \cite{Pit1}, p.471) in
 $u\in[c,d]$, where $[c,d]$ is any finite interval.
 This  implies  the claim (iii) of the theorem.

  Now we turn to the proof of  Theorem  4.6 which will be done  in the  manner similar to the proof of
  Theorem 4.1. For  a given $q\ge 1,$ we consider the stratification points $M_j, j=0,\ldots, q+1$  as in \refm[inc], and
  define the quantities $S^\star_j, E^\star_j,V^\star_j, \ j=0,\ldots, q,$ as in \refm[10] and \refm[19]. Clearly,
  Lemma 5.1 remains valid  with the above points $M_j$. In accordance with the problem considered, the vectors
  $\overrightarrow{L}^{\star}$ and $\overrightarrow{N}^{\star}$ are taken now in a slightly different form: \be
   L_j^{\star}=S_j^{\star}+x_{2j}\sqrt{S_j^{\star}},\quad j=0,\ldots,q
  \la{314}
  \end{equation}
  and
  \be
   N_j^{\star}=E_j^{\star}+x_{2j+1}\sqrt{V_j^{\star}},\quad j=0,\ldots,q-1,\quad
   N_q^{\star}=N-\sum_{j=0}^{q-1}{N_j^{\star}},
  \la{315}
  \end{equation}
  where $x_0,\ldots,x_{2q}$ are arbitrary reals not depending on $N$. Under $\delta=\delta_N$ chosen as in \refm[36], the
 following analog of Lemma 5.3 is then valid:
   \ber
   R(\overrightarrow{L^{\star}},\overrightarrow{N^{\star}})&\sim& \frac{W}{{(2\pi)}^{\frac{q+1}{2}}}
  \left(\prod_{j=0}^{q}{\frac{1}{\sqrt{S_j^{\star}}}}\right)\exp\left({-\half\sum_{j=0}^{q}{x_{2j}^2}}\right),\quad
   N\rightarrow\infty,
   \la{52}
   \ena
   where $W$ is defined by \refm[121]. It is important to note that here, in contrast to  the previous setting,
  \be
  \ds\lim_{N\rightarrow\infty}{\delta M_j}=0,\quad j=0,\ldots,q,
  \la{318}
  \end{equation}
 \nin by \refm[36] and \refm[inc]. Next we will intend to show that the local limit theorems \refm[29] are valid with
 $\kappa_j$ replaced with $\tau_j$  given by \refm[con] below.  Firstly,
 by the definition \refm[inc] of $M_j$ and by
 \refm[318], the integral test gives  for $r\ge 0$ and $N\to\infty,$
  \ber
  & &\sum_{k=M_j}^{M_{j+1}-1}{a_kk^re^{-\delta k}}\sim
  \begin{cases}
    \frac{M_{j+1}^{p+r}-M_j^{p+r}}{p+r}, & \text{if} \ j=0,\ldots,q-1\\
    \Gamma(p+r)\delta^{-p-r}, & \text{if}\ j=q.
  \end{cases}
   \la{317}
  \ena
 \nin From this and \refm[8] we derive that, under $N\to \infty$ and $r>0,$
 \ber
 \E\beta_j^r=
  \begin{cases}
     O\left(M_{j+1}^r\right), & \text{if}\ j=0,\ldots,q-1,\\
     O(\delta^{-r}), & \text{if}\ \ j=q.
  \end{cases}
   \la{324}
   \ena
   In the case considered the  analog of \refm[47] will be
   \ber
    t_{0,j}=
  \begin{cases}
    {\left(L_j^{\star}\right)}^{-\half}M_{j+1}^{-1}\log^2{L_j^{\star}}, & \text{if}\ j=0,\ldots,q-1\\
   \delta{\left(L_q^{\star}\right)}^{-\half}\log^2L_q^{\star}, & \text{if}\ j=q.
  \end{cases}
  \la{325}
  \ena
   Thus, we  obtain that  for any $r\geq0$ and  $N\to\infty,$
   \be
   \left(\E{\beta_j^r}\right)t_{0,j}^r=
   \begin{cases}
    O\left(M_{j+1}^{-\frac{pr}{2}}\log^{2r}{M_{j+1}}\right), & \text{if}\ j=0,\ldots,q-1,\\
    O\left(\delta^{\frac{pr}{2}}\log^{2r}{\delta}\right), & \text{if}\ j=q.
   \end{cases}
   \la{326}
   \end{equation}
   \nin It is easy to see that \refm[326] guarantees the conditions \refm[43]
   and \refm[44].
  Next, we have  from \refm[317]
  \begin{small}
  \ber
  & &\ds{\lim_{N\to\infty}\frac{E_j^{\star}}{\sqrt{V_j^{\star}S_j^{\star}}}}=  \begin{cases}
   \lim_{N\to\infty}\big(\frac{\sqrt{p(p+2)}}{p+1}\big)
   \frac{M_{j+1}^{p+1}-M_j^{p+1}}{\sqrt{\left(M_{j+1}^{p+2}-M_j^{p+2}\right)
   \left(M_{j+1}^{p}-M_j^p\right)}}:= \alpha_j, & \text{if} \ j=0,\ldots,q-1,\\
 \frac{\Gamma(p+1)}{\sqrt{\Gamma(p)\Gamma(p+2)}}:=\alpha_q, & \text{if}\ j=q,
 \end{cases}
 \la{330}
 \ena
 \end{small}
 \nin where the limits $0<\alpha_j< 1,\ j=0,\ldots,q,$ by
 \refm[inc] and some algebra.
  It follows from \refm[315] that
 \ber
 N_q^{\star}=E_q^{\star}-\sum_{j=0}^{q-1}{x_{2j+1}\sqrt{V_j^{\star}}}.
 \la{001}
 \ena
 \nin  We obtain from \refm[317], \refm[inc], \refm[315] and \refm[314]
 \ber
 \non & &\frac{N_j^{\star}-\E{\yn}}{\sqrt{\Var{\yn}}}=\frac{N_j^{\star}-\frac{E_j^{\star}}{S_j^{\star}}L_j^{\star}}
 {\sqrt{\Var{\yn}}}=x_{2j+1}\frac{\sqrt{V_j^{\star}}}{\sqrt{\Var{\yn}}}-x_{2j}\frac{E_j^{\star}}{\sqrt{S_j^{\star}\Var{\yn}}}\\
 & &\non \sim \frac{x_{2j+1}\sqrt{V_j^{\star}}-x_{2j}\frac{E_j^{\star}}{\sqrt{S_j^{\star}}}}
 {\sqrt{V_j^{\star}-\frac{{\left(E_j^{\star}\right)}^2}{S_j^{\star}}}}\to\frac{x_{2j+1}-\alpha_j x_{2j}}
 {\sqrt{1-\alpha_j^2}}:=\tau_j,\\ & &   N\rightarrow\infty,\quad j=0,\ldots,q-1,
 \la{con}
 \ena
 while for $j=q$ we have
 \ber
  \non &
  &\frac{N_q^{\star}-\E{Y^{(q)}}}{\sqrt{\Var{Y^{(q)}}}}=-\frac{\ds\sum_{j=0}^{q-1}{x_{2j+1}\sqrt{V_j^{\star}}}}
  {\sqrt{\Var{Y^{(q)}}}}-x_{2q}\frac{E_q^{\star}}{\sqrt{S_q^{\star}\Var{Y^{(q)}}}}\\
  & & \sim -x_{2q}\frac{E_q^{\star}}{\sqrt{S_q^{\star}\Var{Y^{(q)}}}}\to -\frac{\alpha_q x_{2q}}{\sqrt{1-\alpha_q^2}}
  :=\tau_q,\quad N\to \infty.
 \la{0001}
 \ena
 We  observe that $\vert\tau_j\vert<\infty, \ j=0,\ldots,q.$
   Proceeding further along the same lines as in the proof of Lemma 5.4 leads to the desired analog of
 \refm[29]. Consequently, \refm[065] takes the form
 \begin{small}
 \ber
 R(\overrightarrow{L^{\star}},\overrightarrow{N^{\star}})\sim\frac{\sqrt{\Var{Z_N}}}{{(2\pi)}^{\frac{2q+1}{2}}}
 \left(\prod_{j=0}^{q}{\frac{1}{\sqrt{S_j^{\star}\Var Y^{(j)}}}}\right)
 \exp{\left(-\half\left[\sum_{j=0}^q{x_{2j}^2+{\tau_j^2}}\right]\right)},
 \ N\to \infty.
 \la{116}
 \ena
 \end{small}
    Now we define the scaled quantities for the problem
    considered
    \ber
   \non\hat{\nu}_j^{\star}&=&\frac{\nu_j^{\star}-S_j^{\star}}{\sqrt{S_j^{\star}}},\qquad
    j=0,\ldots,q,\\
\non    \quad
\hat{K}_j^{\star}&=&\frac{K_j^{\star}-E_j^{\star}}{\sqrt{V_j^{\star}}},\qquad
    j=0,\ldots,q-1
    \ena
    and denote
    \ber
\non
\overrightarrow{\hat{\nu}^{\star}}=(\hat{\nu}_0^{\star},\ldots,\hat{\nu}_q^{\star}),\quad
    \overrightarrow{\hat{K}^{\star}}=(\hat{K}_0^{\star},\ldots,\hat{K}_{q-1}^{\star}).
    \la{404}
    \ena

 \nin   \refm[330] and  \refm[62]  imply that
\be \non   \Var{Z_N}\sim  V_q^{\star} \quad \text{and}
\end{equation}
   \be
   \frac{V_j^{\star}}{\Var{Y^{(j)}}}\to \frac{1}{1-\alpha_j^2}<\infty,\quad N\rightarrow\infty,\quad
   j=0,\ldots,q.
   \la{334}
   \end{equation}

    Substituting  \refm[334] into
   \refm[116]  gives \ber
   \non & &\PP\left(\overrightarrow{\hat{\nu}^{\star}}=\overrightarrow{x},\overrightarrow{\hat{K}^{\star}}
   =\overrightarrow{y}\right)\sim{(2\pi)}^{-\frac{2q+1}{2}}
   \left(\prod_{j=0}^{q-1}{\frac{1}{\sqrt{S_j^{\star}V_j^{\star}}}}
   \right)\frac{1}{\sqrt{S_q^{\star}}}\left(\prod_{j=0}^{q}{\frac{1}{\sqrt{1-\alpha_j^2}}}\right)\times\\ & &
    \exp{\left(-\half\left[\sum_{j=0}^{q-1}\left({x_{2j}^2+\frac{\left(x_{2j+1}-\alpha_jx_{2j}\right)^2}{1-\alpha_j^2}}
    \right)+\frac{x_{2q}^2}{1-\alpha_q^2}\right]\right)},\quad
    N\rightarrow\infty,
   \la{336}
   \ena
   where we set $\overrightarrow{x}=(x_0,x_2,\ldots,x_{2q})$ and  $\overrightarrow{y}=(x_1,x_3,\ldots,x_{2q-1})$.\\
\begin{remark}
At the first glance, the expression \refm[336] and its analog
\refm[200] for Theorem 4.1 look very much alike. A crucial
difference between them is that in \refm[200] the variables
$x_j^{\circ}$ are linked via the relation
$\ds\sum_{j=q+1}^{2q+1}{x_j^{\circ}}=0,$ while in \refm[200] the
variables $x_j$ are free. Formally, this is implied by the fact
that the value $\kappa_q$ in \refm[50] depends on $x_{q}^{\circ}$
and $x_{2q+1}^{\circ},$ whereas the value $\tau_q$ in \refm[0001]
depends on $x_{2q}$ only. As a result of the aforementioned
difference, conditional independence in Theorem 4.1 transforms to
independence in Theorem 4.3.
\end{remark}
 \nin From \refm[336], the claim (i) of Theorem 4.6 follows immediately.
 It is left to find  explicitly the covariance matrix of the corresponding limiting gaussian
  distribution.
  For the sake of convenience, we write
 $\vec{z}=(x_0,\ldots,x_{2q}).$ Then the matrix, say $B$, of the quadratic form in the exponent of
 \refm[336] has  a diagonal block structure with the blocks $B_j$ of the form
 \ber
\non B_j=\left( \begin{array}{cccc}
     \frac{1}{1-\alpha_j^2}  &-\frac{\alpha_j}{1-\alpha_j^2}  & \\
     -\frac{\alpha_j}{1-\alpha_j^2} & \frac{1}{1-\alpha_j^2}
     \end{array} \right),\quad j=0,\ldots,q-1
 \la{342}
 \ena
 and
 \be
  B_q= \frac{1}{1-\alpha_q^2}. \la{143}
 \end{equation}
  Thus,
  \be
   \ds\det{B}=\prod_{j=0}^{q}{\det{B_j}}=\prod_{j=0}^{q}{\frac{1}{1-\alpha_j^2}}
   \la{338}
   \end{equation}
 \nin and
 \be
  B_j^{-1}=\left( \begin{array}{cccc}
     1 &\alpha_j  & \\
     \alpha_j & 1
     \end{array} \right),\quad j=0,\ldots,q-1,\qquad B_q^{-1}=1-\alpha_q^2.
    \la{cov1_new}
    \end{equation}
 \nin In accordance with the above notation for $\vec{z}$, we set                                     $\overrightarrow{\hat{\chi}}=(\hat{\nu}_0^{\star},\hat{K}_0^{\star},      \hat{\nu}_1^{\star},\hat{K}_1^{\star},\ldots,\hat{\nu}_{q-1}^{\star},      K_{q-1}^{\star},\hat{\nu}_q^{\star})$ which
   enables to rewrite \refm[336] as
 \ber
 \non & &\PP\left(\overrightarrow{\hat{\chi}}=\vec{z}
 \right)\sim\frac{1}{{(2\pi)}^{\frac{2q+1}{2}}\sqrt{\det{B}}}\left(\prod_{j=0}^{q-1}
 {\frac{1}{\sqrt{S_j^{\star}V_j^{\star}}}}\right)\frac{1}{\sqrt{S_q^{\star}}}\exp{\left(-\half\overrightarrow{z}B
 {\overrightarrow{z}}^T\right)},\\ & &N\rightarrow\infty.
 \la{341}
 \ena
 Hence, the covariance matrix $B^{-1}$ of the gaussian distribution can be written
 as the  Kronecker product of blocks $B_j^{-1}$ as in \refm[cov1_new]:

 \ber
 B^{-1}=B_0^{-1}\otimes B_1^{-1}\otimes\ldots\otimes B_q^{-1}.
 \la{cov11}
 \ena
 \nin Finally, by the same argument as before, we derive the weak convergence as stated in {\rm(ii)}   of Theorem 4.6.
  \section{ Comparison with limit shapes for the generalized Bose-Einstein
  and Fermi-Dirac models of ideal gas}

Recall (see Section 3) that the generalized BE and the FD models
  in the title are
  associated with
  the  generating functions $\nicef^{(BE)}, \nicef^{(FD)}$ respectively:
 \ber
\non \nicef^{(BE)}(x)=\prod_{k\ge 1}\frac{1}{(1-x^k)^{m_k}},\quad
\nicef^{(FD)}(x)=\prod_{k\ge 1}(1+x^k)^{m_k},\\ \quad \vert
x\vert<1,\quad m_k=Ck^{p-1},
 \quad p\ge 1,\quad C>0.
 \la{090}
 \ena
Limit shapes for these models, say $C_{p-1}^{(BE)},C_{p-1}
^{(FD)}$ respectively, were obtained by Vershik in
\cite{V1},\cite{V23}: \be\la{lshc}
C_{p-1}^{(BE)}(u)=\int_{u}^\infty
x^{p-1}\frac{e^{-c_1x}}{1-e^{-c_1x}}dx,\quad
C_{p-1}^{(FD)}(u)=\int_{u}^\infty
x^{p-1}\frac{e^{-c_2x}}{1+e^{-c_2x}}dx,\quad u\ge 0,\ p\ge 1,
\end{equation} where $c_1=c_1(p)>0,\ c_2=c_2(p)>0$ are  normalizing constants.

 It
is interesting to compare these limit shapes with the ones,
denoted $l_{p-1},$
  in our Theorem 4.3:

\be l_{p-1}(u)=\frac{\int_u^\infty{e^{-x}x^{p-1}dx}}{\Gamma(p+1)},
\quad u\ge 0,\quad p>0.\la{lshc1}
\end{equation}

\nin  Firstly,   note that the limit shapes for  both models
  are derived under the   scalings $r_N$ of the same order $N^{\frac{1}{p+1}}.$
  In view of Proposition 2.5, the immediate explanation of this striking
        coincidence is that by the result of Vershik and Yakubovich (\cite{VY1},
         Section 5), the models \refm[090] of ideal gas have the same threshold $N^{\frac{1}{p+1}}$
         as
         the models
considered in the paper.
         In a broader context, we point
that the generalized BE models and the models \refm[90] studied in
our paper, are linked with each other via the exponentiation of
the generating function $\nicef^{(BE)}$ (see \cite{Bur}).
Moreover, the Bell-Burris Lemma 5.1 in \cite{BB} states that if
the parameter function $m=\{m_k\}$ for the generalized BE model
(multiset) is such that $\lim_{k\to\infty}\frac{m_{k-1}}{m_k}=y$
with $y<1,$
  then the corresponding multiplicative measure for the  BE model is  asymptotically equivalent to the
 measure $\mu_N$ (assembly) given by \refm[90] with $a_k\sim m_k, \ k\to \infty.$
   However, for the generalized BE models \refm[090],
   $\lim_{k\to\infty}\frac{m_{k-1}}{m_k}=1,$
     which explains the following basic
     difference in the form of the limit curves $C_{p-1}^{(BE)}$ and $l_{p-1}$.
 From \refm[lshc1],\refm[lshc] we see  that \be
l_{p-1}(0)
  =p^{-1}<\infty,\quad p>0,\la{001_new}\eu while
\be C_{0}^{(BE)}(0)=\lim_{u\to 0^+} -\frac{\sqrt{6}}{\pi}
 \log{\left(1-\exp\left({-\frac{\pi u}{\sqrt{6}}}\right)\right)}=\infty,\quad C_{p-1}(0)<\infty, \quad p>1.
\la{0907}\eu
 By \refm[85],\refm[86] and \refm[1],
   the value of a limit shape at $u=0$ "approximates"
   the random variable $\widetilde{\nu}(0)=$ the total number of components
in a random partition,      multiplied by the factor
$\frac{r_N}{N}.$ In  \cite{V1},\cite{V23}   the phenomenon
\refm[0907] is linked   with the Bose -Einstein condensation of
energy. According to this interpretation, the finiteness of the
limit shape at $u=0$ indicates the appearance of
 condensation of energy (around  the value $\frac{N}{r_N}$),
 and in view of this the value $p=1$ (=uniform distribution on the set $\Omega_N$)
 was distinguished  as the phase
transition point for the generalized BE models \refm[090]. We now
give an analytic explanation of  the fact that
$C_{0}^{(BE)}(0)=\infty$ in contrast to $l_0(0)<\infty.$ By the
seminal result of Erd$\ddot{o}s$ and Lehner(1941) (for references
see \cite{VY1}), the number of components in a random partition of
$N$ is asymptotically $\frac{2\pi}{\sqrt{6}}\sqrt{N}\log{N}\gg
\sqrt{N}=\frac{N}{r_N}.$ On the other hand, it was proven in
\cite{EG} that the number of components in our model with $p=1$ is
asymptotically $\sqrt{N}=h\frac{N}{r_N},$ where the constant $h$
is as in \refm[h].

 \vskip 2cm
 \nin {\bf
Acknowledgement}

 We appreciate the illuminating remarks and the
constructive suggestions of a referee.

  
\end{document}